\documentclass[11pt,a4paper]{amsart}

\title[Koppelman formulas on Grassmannians]{Koppelman formulas on Grassmannians}

\author{Elin G\"{o}tmark \& H\aa kan Samuelsson \& Henrik Sepp\"{a}nen}
\thanks{The second author was supported by a Post Doctoral Fellowship from the Swedish Research Council.}


\address{E. G\"{o}tmark, Department of Mathematical Sciences, Division of Mathematics, Chalmers University of Technology
and G\"{o}teborg University, SE-412 96 G\"{o}teborg, Sweden}
\email{elin@math.chalmers.se}
\address{H. Samuelsson, Department of Mathematics, University of Wuppertal, Gaussstrasse 20, D-42119 Wuppertal, Germany}
\email{hasam@math.chalmers.se}
\address{H. Sepp\"{a}nen, Department of Mathematical Sciences, Division of Mathematics, Chalmers University of Technology
and G\"{o}teborg University, SE-412 96 G\"{o}teborg, Sweden}
\email{henriks@math.chalmers.se}

\usepackage[english]{babel}
\usepackage{amsmath}
\usepackage{amsthm,amssymb,latexsym}
\usepackage{t1enc}
\usepackage{mathrsfs}
\usepackage{graphicx}
\usepackage{xypic}

\newtheorem{proposition}{Proposition}
\newtheorem{theorem}[proposition]{Theorem}
\newtheorem{lemma}[proposition]{Lemma}
\newtheorem{corollary}[proposition]{Corollary}

\theoremstyle{definition}
\newtheorem{definition}[proposition]{Definition}

\newtheorem{remark}[proposition]{Remark}

\newcommand{\C}{\mathbb{C}}

\newcommand{\debar}{\bar{\partial}}

\newcommand{\PV}{\mathcal{P}}

\newcommand{\CP}{\mathbb{CP}}

\newcommand{\z}{\zeta}
\newcommand{\deta}{\delta_{\eta}}

\newcommand{\Hom}{\textrm{Hom}}
\newcommand{\Hram}{\mathfrak{h}}
\newcommand{\Fram}{\mathfrak{f}}
\newcommand{\Eram}{\mathfrak{e}}

\def\newop#1{\expandafter\def\csname #1\endcsname{\mathop{\rm #1}\nolimits}}
\newop{span}

\begin{document}
\nocite{*}
\bibliographystyle{plain}

\begin{abstract}
We construct Koppelman formulas on Grassmannians for forms with values 
in any holomorphic line bundle as well as in the tautological vector bundle and 
its dual.
As a consequence we obtain some vanishing theorems of the Bott-Borel-Weil type.
We also relate the projection part of our formulas to 
the Bergman kernels associated to the line bundles.
\end{abstract}

\maketitle
\thispagestyle{empty}

\section{Introduction}

The Cauchy integral formula in one complex variable is
of vast
importance in many respects. It provides a way of representing a
holomorphic function as a superposition of simple rational functions,
and gives an explicit solution to the equation $\debar u =
f$. Furthermore, it is an important tool in function theory.
For our purposes it is convenient to
note that Cauchy's formula is equivalent to the current equation
$\debar u=[z]$, where $u=(2\pi i)^{-1}d\z/(\z-z)$ is the Cauchy form, and $[z]$ is the Dirac measure at 
$z$ considered as a $(1,1)$-current. This point of view is well adapted 
for generating weighted Cauchy formulas. For instance, by computing 
$\debar \big(((1-|\z|^2)/(1-z\bar{\z}))^{\alpha} u\big)$ 
in the current sense, one obtains (for suitable $\alpha$) 
the weighted representation 
formula

\begin{equation*}
f(z)=\frac{\alpha}{\pi}\int_{\{|\z|<1\}}
f(\z)\frac{(1-|\z|^2)^{\alpha-1}}{(1-z\bar{\z})^{\alpha+1}}d\lambda(\z),
\end{equation*}
for holomorphic functions on the unit disc with certain limited growth at the boundary.
The integral kernel is the reproducing kernel for a weighted Bergman space; and this
shows that there is a connection between Cauchy kernels and Bergman kernels. Both these kernels
are also intimately linked with the symmetry of the disc.
Recall that the group 

\begin{eqnarray*}
SU(1,1)
=\left\{\left( \begin{array}{cc}
a & b\\
\overline{b} & \overline{a} \end{array}\right) \in M_{22}(\mathbb{C}) |\; |a|^2-|b|^2=1\right\}
\end{eqnarray*}
acts holomorphically and transitively on the unit disc by $z \mapsto (az+b)/(\overline{b}z+\overline{a})$.
The stabilizer of the origin is the subgroup

\begin{equation*}
K:= \left\{\left( \begin{array}{cc}
e^{i\theta} & 0\\
0 & e^{-i\theta} \end{array}\right)\right\}\cong S^1,
\end{equation*} 
and hence the disc can be viewed as the homogeneous space
$SU(1,1)/S^1$.  
The kernels are then invariant under certain actions on functions which are induced from the 
natural action on the closed disc.
From the point of view of representation theory, 
the Bergman kernels are interesting since the corresponding weighted Bergman spaces form a family
of unitary representation spaces for $SU(1,1)$, and moreover, these kernels can be described entirely in terms of
the Lie-theoretic structure of the group.
This discussion indicates two possible directions of generalizations; namely to
domains in $\C^n$, and to complex homogeneous spaces. In the latter case, the class of bounded symmetric domains
have been studied extensively from the Lie-theoretic point of view. Hua, \cite{hua}, computed the Cauchy kernels
and Bergman kernels for the classical domains using the explicit description of their symmetry groups.
Later, more abstract group theoretic machinery has been used to describe both Bergman kernels (cf.\ \cite{satake}) and
the generalized Cauchy-Szeg kernels, \cite{K-W}. For compact Hermitian symmetric spaces, Bergman kernels for line bundles
can be described explicitly in terms of the polynomial models for the spaces of global holomorphic sections, \cite{genkai}.

Complex analysts have mainly been concerned with domains in $\C^n$. The Bochner-Martinelli kernel
represents holomorphic functions in any domain but has the drawback of
not being holomorphic, a property which is highly useful in
applications. The Cauchy-Fantappi-Leray
kernel is holomorphic in domains where we can find a holomorphic
support function, for example strictly pseudoconvex domains.  
More flexibility is afforded by using weighted formulas, which was
first done in \cite{AB}, and such formulas have been widely used in
applications such as interpolation, division, obtaining estimates for
solutions to the $\debar$-equation, etc. See, e.g., \cite{MA1} and
\cite{BY2} and the references therein. 
Some work has also been done on generalizing integral formulas to complex manifolds,
see, e.g., \cite{HL}, \cite{BE2}, \cite{bob1}. Of these, the paper \cite{bob1} by 
Berndtsson will be of particular importance for us; see below.

More recently, in \cite{MA1} was introduced a general method for
generating weighted formulas for domains in $\mathbb{C}^n$, both for holomorphic
functions and $(p,q)$-forms. For future reference,
we will describe this method in the former case in some detail. First, recall
that the Cauchy kernel, $u$, in 
one variable satisfies $\debar u = [z]$, but less obviously, we
also have $\delta_{\z-z} u = 1$, 
where $\delta_{\z-z}$ denotes contraction with the vector field $2 \pi i
(\z-z) \partial /\partial \z$. These equations can be combined into the
single equation 

\begin{equation} \label{zira}
\nabla_{\z-z} u = 1-[z],
\end{equation}
where $\nabla_{\z-z}$ is the operator 

\begin{equation*}
\nabla_{\z-z}=\delta_{\z-z}-\debar.
\end{equation*}
To generalize this to $\mathbb{C}^n$, we define $\delta_{\z-z}$ as
contraction with 
\begin{equation} \label{seseli}
2 \pi i\sum (\z_j-z_j) \frac{\partial}{\partial \z_j},
\end{equation}
and if we construe equation \eqref{zira} as being in $\mathbb{C}^n$, the right
hand side of \eqref{zira} now contains one form of 
bidegree $(0,0)$ and one of bidegree $(n,n)$, so we must in fact have $u
= u_{1,0} + u_{2,1} + \ldots + u_{n,n-1}$, where $u_{k,k-1}$ has
bidegree $(k,k-1)$. We can then write the $\nabla_{\z-z}$-equation (\ref{zira}) as the system 
of equations 

\begin{equation*}
\delta_{\z-z} u_{1,0} = 1, \qquad \delta_{\z-z} u_{1,2} - \debar u_{1,0} = 0, \qquad 
\ldots, \qquad \debar u_{n,n-1} = [z]. 
\end{equation*}
In that case, $u_{n,n-1}$ will
satisfy $\debar u_{n,n-1} = [z]$ and will give a kernel for a
representation formula.  
One advantage of this approach, as opposed to just solving $\debar
u_{n,n-1} = [z]$, is that it easily allows for weighted integral
formulas. We define $g = g_{0,0} + \cdots + g_{n,n}$ to be a weight if
$\nabla g = 0$ and $g_{0,0}(z,z) = 1$. It is easy to see that $\nabla (u \wedge
g) = g - [\Delta]$, and this yields a representation formula

\begin{equation*}
\phi(z) = \int_{\partial D} \phi(\z) (u \wedge g)_n + \int_D \phi g_n
\end{equation*}
if $\phi \in \mathcal{O}(\overline D)$ and $z \in D$. Note that
if $g_1$ and $g_2$ are weights, then $g_1 \wedge g_2$ is also a weight. 

In the case of compact manifolds one is naturally led to consider holomorphic line bundles
and representation formulas for holomorphic sections as well as smooth bundle-valued forms.
In this setting the integral kernels must be operator valued, and the 
integrals become superpositions of contributions from all fibres. Our method for achieving this 
has two crucial components; the above mentioned $\nabla$-formalism, and Berndtsson's method
from \cite{bob1}. Indeed, Berndtsson gave a method for obtaining integral 
formulas for $(p,q)$-forms on $n$-dimensional manifolds $X$ which admit a vector bundle
of rank $n$ over $X \times X$ such that the diagonal has a defining
section $\eta$; and to get formulas for forms with values in bundles the $\nabla$-method 
is well suited. In fact, by generalizing it to manifolds one realizes that it allows for 
operator valued weights. We then need 
something to 
substitute for the vector field \eqref{seseli}, and this is where Berndtsson's
assumption comes in: we will use the section $\eta$ to contract with,
and define $\nabla_{\eta} := \delta_{\eta} - \debar$. 
It is of independent interest
to note that $\nabla_{\eta}$ in fact is a superconnection in the sense of Quillen, \cite{QU}.
In the recent article \cite{elin} by the first author, this general theory for integral formulas on manifolds
has been developed to a large degree,
and explicit formulas have been constructed on $\CP^n$ yielding explicit proofs of vanishing theorems for 
its line bundles. Such proofs could be of interest also for representation theoretic purposes. Indeed, 
in view of the by now firmly established
goal, initiated by the Bott-Borel-Weil theorem and further fortified by the conjecture of Langlands, 
\cite{lang}, and
Schmid's proof of it, \cite{schmid-L^2}, of wanting to realize representations of Lie groups in 
Dolbeault cohomology (or, rather
$L^2$-cohomology in the non-compact case), (cf.\ also \cite{wong1} and \cite{wong2}), 
it is our hope that explicit integral formulas could give further insight into the underlying group theory.     

In this paper, we extend the method in \cite{elin} to the vector bundle setting and we apply the technique 
to complex Grassmannians, $Gr(k,N)$.
We find a suitable vector bundle, with a section $\eta$ as above, and natural weights for the line bundles and for the 
tautological $k$-plane bundle.
We thus get Koppelman formulas for $(p,q)$-forms with values in any holomorphic line bundle as well as in 
the tautological bundle and its dual.
The construction is uniform in the sense that it uses the explicit description of the Picard group
of holomorphic line bundles and reduces the problem to that of finding a weight for the generator.
The generator in turn, is the determinant of the tautological bundle; by certain algebraic properties of weights,
it thus suffices to construct a weight for the tautological bundle.  
As an application, we give explicit proofs of certain vanishing theorems of Bott-Borel-Weil type
\footnote{These are not given in the form including the $\rho$-shift which is common in representation theory.}  
for the cohomology groups associated with these line bundles. We also
relate the projection part of our Koppelman formulas to Bergman kernels; thus giving a geometric interpretation 
of the latter ones.

This paper is organized as follows: In Section \ref{yuuri} we recapture the general method for finding weighted 
Koppelman formulas on manifolds from \cite{elin}. The only difference is that we allow
for forms with values in vector bundles and state a slightly more general Koppelman formula. The proofs 
have been omitted since they are straightforward generalizations of the proofs in \cite{elin}. Section 
\ref{viktsektion} describes some general operations on weights. In Section \ref{yuuuri} we construct the 
ingredients necessary
to generate weighted formulas on Grassmannians according to the general framework. 
In Section \ref{yuuuuri} we review the representation theoretic description of the Picard group and we prove a certain
invariance property for the weights, which will be useful for the applications. We also prove that 
the bundle $E$ restricted to the diagonal is equivalent to the holomorphic cotangent bundle over $Gr(k,N)$. 
In the last section, Section \ref{yuuuuuri}, we discuss some applications; we obtain vanishing theorems 
for the line bundles over Grassmann, and we give a geometric interpretation of the Bergman kernels associated to 
the line bundles.

{\bf Acknowledgement:}
We are grateful to Mats Andersson and Genkai Zhang for rewarding
discussions and for valuable comments on preliminary versions of this
paper. We would also like to thank Harald Upmeier for interesting
discussions on the topic of this paper.


\section{A general method for finding weighted Koppelman formulas on manifolds}\label{yuuri}

Let $X$ be a complex manifold of dimension $n$. We want to find
Koppelman formulas for differential forms on $X$ with values in a given vector
bundle $H \rightarrow X$. The method described in this section is 
taken from \cite{elin}, except for the generalization which yields
formulas for a general vector bundle $H$ instead of for a line
bundle. 

We begin by noting that Stokes' theorem holds also for
sections of vector bundles, which is easily proved. 
Let $M$ be any complex manifold, and $G \rightarrow M$ a holomorphic Hermitian
vector bundle over $M$. Let $D_{G^{\ast}}$ and $D_G$ be the Chern connections
for $G^{\ast}$ and $G$ respectively. If $u$ is a differential form taking values
in $G^{\ast}$ and $\phi$ is a test form with values in $G$, we have 

\begin{equation}\label{stenbrott}
\int_M D_{G^{\ast}} u \wedge \phi = (-1)^{\deg u+1} \int_M u \wedge D_G \phi,
\end{equation}
where $\wedge$ denotes taking the natural pairing between the factors
in $G^{\ast}$ and $G$, and taking the wedge product between the
factors which are differential forms.
If $u$ is instead a current, we can take \eqref{stenbrott} as a
definition. In the same way, we also have

\begin{equation}\label{havtorn}
\int_M \debar u \wedge \phi = (-1)^{\deg u+1} \int_M u \wedge \debar \phi.
\end{equation}

Let $\Delta$ be the diagonal in $X_{z} \times X_{\z}$. Let $H_{z}$ denote
$\pi_{z}^{*}(H)$, where $\pi_{z}$ 
is the projection from $X_{z} \times X_{\z}$ to $X_{z}$, and analogously
for $H_{\z}$. Let $g_{0,0}$ be a section of $H_z \otimes H_{\z}^{\ast} =
\textrm{Hom}(H_{\z},H_z)$ such that $g_{0,0}(z,z) = \textrm{Id}$ for all
$z$. If $[\Delta]$ denotes the current of integration over the
diagonal and $\omega(\z,z)$ is a 
differential form with values in $H_z^{\ast} \otimes H_{\z}$, then we let

\begin{equation*}
[\Delta]_{g_{0,0}} (\omega): = [\Delta]. ((g_{0,0}
\otimes \textrm{Id}) \omega),
\end{equation*}
where $\textrm{Id}$ acts on the differential forms in $\omega$, and we
take the natural pairing $(H_z^{\ast} \otimes H_{\z}) \times (H_z
\otimes H_{\z}^{\ast}) \to \mathbb{C}$. Note
that this does not depend on which $g_{0,0}$ we choose, since the
values on the diagonal are the only ones that matter. The reason for
the subscript on $g_{0,0}$ will become apparent later on. 

\begin{proposition}[Koppelman's formula]\label{koala}
Assume that $D \subset X_{\z}$, $\phi \in \mathcal{E}_{p,q}
(\bar{D}, H_{\z})$, and that the current $K(z,\z)$ and the smooth form
$P(z,\z)$ take values in $H_z \otimes H_{\z}^{\ast} = \textrm{Hom}(H_{\z},H_z)$ and
solve the equation 
\begin{equation}\label{miyazaki}
\debar K = [\Delta]_{g_{0,0}} - P.
\end{equation}
We then have 
\begin{equation} \label{iris}
\phi(z) =  \int_{\partial D} K \wedge \phi + \int_{D} K \wedge \debar \phi +  
\debar_z \int_{D} K \wedge \phi + \int_{D} P \wedge \phi,
\end{equation}
where the integrals are taken over the $\z$ variable. 
\end{proposition}

The proof of this uses \eqref{havtorn} but is otherwise
just like the usual proof of the Koppelman formula. Note that if
$\phi$ in \eqref{iris} is a $\debar$-closed form and the  
first and fourth terms of the right hand side 
of Koppelman's formula vanish, we get a 
solution to the $\debar$-problem for $\phi$. 

Our purpose now is to find $K$ and $P$ that satisfy \eqref{miyazaki}
in a special type of manifold. To begin with, we will let $H$ be the
trivial line bundle. Assume that we can find a holomorphic vector
bundle $E \rightarrow 
X_{z} \times X_{\z}$ of rank $n$, such that there exists 
a holomorphic section $\eta$ of $E$ that 
defines the diagonal $\Delta$.
In other words, $\eta$ must vanish to the first order on $\Delta$ and be
non-zero elsewhere. Let $\{e_i\}$ be a local frame for $E$, and
$\{e_i^{\ast}\}$ the dual local frame for $E^{\ast}$. 
Contraction with $\eta$ is an 
operation on $E^{\ast}$ which we denote by $\deta$; if $\eta = 
\sum \eta_i e_i$ then

\begin{equation*}
\deta \left(\sum \sigma_i e_i^{\ast} \right) = \sum \eta_i \sigma_i.
\end{equation*}
We define the operator 
\begin{equation*}
\nabla_{\eta} = \deta - \debar.
\end{equation*}

Choose a Hermitian metric $h$ for $E$, 
let $D_E$ be the Chern connection on $E$, and $D_{E^{\ast}}$ the induced
connection on $E^{\ast}$. Consider the bundle 

\begin{equation*}
G_E =
\Lambda [T^{\ast}(X \times X) \oplus E \oplus
E^{\ast}] \to X \times X
\end{equation*}
and $\Gamma(X \times X, G_E)$, the space
of $C^{\infty}$ sections of $G_E$ (note the change of notation compared
to \cite{elin}). If $A$ lies in $\Gamma (X \times X,
T^{\ast}(X \times X) 
\otimes E \otimes E^{\ast}))$, then we define 
$\tilde{A}$ as the corresponding element in $\Gamma(X \times X, G_E)$,
arranged with the 
differential form first, then the section of $E$ and finally the
section of $E^{\ast}$. For example, if $A =
dz_1 \otimes e_1 \otimes e_1^{\ast}$, then $\tilde{A} = dz_1 \wedge e_1
\wedge e_1^{\ast}$.  

To define a derivation $D$ on $\Gamma(X \times X, G_E)$, we first let $Df
= \widetilde{D_E f}$ for a section $f$ of $E$, and $Dg = \widetilde{D_{E^{\ast}}g}$
for a section $g$ of $E^{\ast}$. We then extend the definition by 

\begin{equation*}
D(\xi_1 \wedge \xi_2) = D  \xi_1 \wedge \xi_2 + (-1)^{\deg \xi_1} \xi_1 \wedge D \xi_2,
\end{equation*}
where $D\xi_i = d\xi_i$ if $\xi_i$ happens to be a differential form,
and $\deg\xi_1$ is the total degree of $\xi_1$. For example,
$\deg(\alpha \wedge e_1 \wedge e_1^{\ast}) = \deg\alpha + 2$, where $\deg\alpha$
is the degree of $\alpha$ as a differential form. 
We let 

\begin{equation*}
\mathcal{L}^m = \bigoplus_p \Gamma (X \times X, \Lambda^p E^{\ast} \wedge
\Lambda^{p+m} T^{\ast}_{0,1} (X \times X));
\end{equation*}
note that $\mathcal{L}^m$ is a subspace of $\Gamma(X \times X, G_E)$. The
operator $\nabla_{\eta}$ will act in a natural way as $\nabla_{\eta} \colon \mathcal{L}^m
\to \mathcal{L}^{m+1}$.  If $f \in
\mathcal{L}^m$ and $g \in \mathcal{L}^{k}$, then $f \wedge g \in
\mathcal{L}^{m+k}$. We also see that $\nabla_{\eta}$ obeys Leibniz' rule, and
that $\nabla_{\eta}^2 = 0$.

\begin{definition} \label{E-integral}
For a form $f(z,\z)$ on $X \times X$, we define

\begin{equation*}
\int_E f(z,\z) \wedge e_1 \wedge e_1^{\ast} \wedge \ldots \wedge e_n \wedge e_n^{\ast} = f(z,\z).
\end{equation*}
\end{definition}
\noindent
Note that if $I$ is the identity on $E$, then $\tilde{I} = e \wedge
e^{\ast} = e_1 \wedge e_1^{\ast} + \ldots + e_n \wedge e_n^{\ast}$. It follows that
$\tilde{I}_n = e_1 \wedge e_1^{\ast} \wedge \ldots \wedge e_n \wedge e_n^{\ast}$ (with the
notation $a_n = a^n/n!$), so the
definition above is independent of the choice of frame. Our derivation
$D$ and $\int_E$ interact in the following way:

\begin{proposition}\label{ek}
If $F \in \Gamma(X \times X, G_E)$ then 
\begin{equation*}
d \int_E F = \int_E DF.
\end{equation*}
\end{proposition}

We will now construct integral formulas on $X \times X$. As a first
step, we find a section $\sigma$ of $E^{\ast}$ such that $\deta 
\sigma = 1$ outside $\Delta$. For reasons that will become apparent,
we choose $\sigma$ to have minimal pointwise norm with respect to the metric
$h$, which means that $\sigma = \sum_{ij} h_{ij} \bar{\eta}_j e^{\ast}_i /|\eta|^2$. 
Close to $\Delta$, it is obvious that
$|\sigma| \lesssim 1/|\eta|$, and a calculation shows that we also
have $|\debar \sigma| \lesssim 1/|\eta|^2$. Next, we construct a
section $u$ with the property that $\nabla_{\eta} u = 1 
- R$ where $R$ is a current with support on $\Delta$. We set 

\begin{equation}\label{nord}
u = \frac{\sigma}{\nabla_{\eta} \sigma} = \sum_{k=0}^\infty \sigma 
\wedge (\debar \sigma)^k,
\end{equation}
and note that $u \in \mathcal{L}^{-1}$. By $u_{k,k-1}$ we will mean
the term in $u$ of degree $k$ in $E^{\ast}$ and degree $k-1$ in
$T^{\ast}_{0,1} (X \times X)$. It is easily checked that $\nabla_{\eta}
u = 1$ outside $\Delta$.

The following theorem yields a Koppelman formula by Theorem
\ref{koala}, with the trivial line bundle as $H$:
\begin{theorem}\label{godis}
Let $E \to X \times X$ be a vector bundle with a section $\eta$ which
defines the diagonal $\Delta$ of $X \times X$. 
We have 

\begin{equation*}
\debar K = [\Delta] - P,
\end{equation*}
where 

\begin{equation}\label{calopteryx}
K = \int_E u \wedge \left(\frac{D \eta}{2 \pi i} + \frac{i
    \tilde{\Theta}}{2 \pi}\right)_n 
\quad \textrm{and} \quad P = \int_E \left(\frac{D \eta}{2 \pi i} +
  \frac{i\tilde{\Theta}}{2 \pi } \right)_n, 
\end{equation}
and $u$ is defined by \eqref{nord}. 
\end{theorem}

\noindent
Note that since $D \eta$ contains no $e_i^{\ast}$'s, we have 
\begin{equation*}
P = \int_E (\frac{i\tilde{\Theta}}{2 \pi })_n = \det\frac{i\Theta}{2
  \pi} = c_n(E), 
\end{equation*}
i.e., the $n$th Chern class of $E$. The factor 

\begin{equation*}
\frac{D \eta}{2 \pi i} + \frac{i\tilde{\Theta}}{2 \pi}
\end{equation*}
is actually the supercurvature associated with the operator $\nabla_{\eta}$ if we
view $\nabla_{\eta}$ as a superconnection in the sense of
Quillen, \cite{QU}. In fact, we have the following Bianchi identity:

\begin{equation} \label{fredag}
\nabla_{\eta} \left(\frac{D \eta}{2 \pi i} + \frac{i\tilde{\Theta}}{2
    \pi}\right) = 0,
\end{equation}
for a direct proof see, e.g., \cite{elin}.

The idea behind the proof of Theorem \ref{godis} is that by
\eqref{fredag} and Proposition \ref{ek} we have 

\begin{eqnarray} \label{angel}
& & \debar \int_E u \wedge \left(\frac{D \eta}{2 \pi i} + \frac{i
    \tilde{\Theta}}{2 \pi}\right)_n = 
\int_E \debar \left[ u \wedge \left(\frac{D \eta}{2 \pi i} + \frac{i
    \tilde{\Theta}}{2 \pi} \right)_n \right] = \nonumber \\
& = & - \int_E \nabla_\eta \left[ u \wedge \left(\frac{D \eta}{2 \pi i} + \frac{i
    \tilde{\Theta}}{2 \pi} \right)_n \right] = \nonumber \\
& = & - \int_E \left(\frac{D \eta}{2 \pi i} + \frac{i
    \tilde{\Theta}}{2 \pi}\right)_n + \frac{1}{(2 \pi i)^n} \int_E R 
\wedge (D \eta)_n. 
\end{eqnarray}
The left hand term in \eqref{angel} is $P$. The rest of the proof
consists of proving that 

\begin{equation}\label{kajflock}
\frac{1}{(2 \pi i)^n} \int_E R \wedge (D \eta)_n = [ \Delta ],
\end{equation}
which is proved by choosing local coordinates
on $X$, and reducing the problem to the $\mathbb{C}^n$-case. For
details of the proof, see, e.g., \cite{elin}. 

As explained in the introduction, we will obtain more flexible
formulas if we use weights. 

\begin{definition}
A section $g$ with values in $\mathcal{L}_0$ is a weight if 
$\nabla_{\eta} g = 0$ and $g_{0,0}(z,z) = 1$.
\end{definition}

\noindent
Theorem \ref{godis} goes through with essentially the same proof if we
take 

\begin{equation}\label{thai}
K_g = \int_E u \wedge g \wedge \left(\frac{D \eta}{2 \pi i} + \frac{i
    \tilde{\Theta}}{2 \pi}\right)_n
\quad \textrm{and} \quad P_g = \int_E g \wedge \left(\frac{D \eta}{2 \pi
    i} + \frac{i \tilde{\Theta}}{2 \pi}\right)_n, 
\end{equation}
as shown by the following calculation: 

\begin{eqnarray}\label{mufti}
\debar K_g & = & -\int_E \nabla_{\eta} u \wedge g \wedge \left(\frac{D \eta}{2 \pi i} + \frac{i
    \tilde{\Theta}}{2 \pi}\right)_n = \nonumber \\
& = & -\int_E (g - R) \wedge \left(\frac{D \eta}{2 \pi i} + \frac{i
    \tilde{\Theta}}{2 \pi}\right)_n = [\Delta ] - P_g, 
\end{eqnarray}
which follows from the proof of Theorem \ref{godis} and the properties of
weights.

Finally, we will use weights taking values in $\textrm{Hom}(H_{\z},H_z)$
to construct Koppelman formulas for differential forms with
values in the vector bundle $H \to X$. We define 

\begin{equation*}
G_{E,H} = \textrm{Hom}(H_{\z},H_z) \otimes \Lambda [T^{\ast}(X \times
X) \oplus E \oplus E^{\ast}] \to X \times X
\end{equation*}
and  
\begin{equation}\label{vikter}
\mathcal{L}^m_H := \bigoplus_p \Gamma (X \times X,
\textrm{Hom}(H_{\z},H_z) \otimes [\Lambda^p E^{\ast} \wedge
\Lambda^{p+m} T^{\ast}_{0,1} (X \times X)]).
\end{equation}
We define $\deta$ on $\Gamma (X \times X, G_{E,H})$ as $\textrm{Id} \otimes
\deta$, where $\textrm{Id}$ acts on the factors in
$\textrm{Hom}(H_{\z},H_z)$ and $\deta$ on the factors in $\Lambda
[T^{\ast}(X \times X) \oplus E \oplus E^{\ast}]$. 
We also need to extend the derivation $D$ to $\Gamma (X \times X,
G_{E,H})$. If $a_1$ is a
differential form taking values in 
$\textrm{Hom}(H_{\z},H_z)$, and $a_2 \in \Gamma (X \times X, G_{E})$,
then we define  

\begin{equation*}
D (a_1 \wedge a_2) = D_{\textrm{Hom}(H_{\z},H_z)} a_1 \wedge a_2 +
(-1)^{\deg a_1} a_1 \wedge D a_2, 
\end{equation*}
where $D_{\textrm{Hom}(H_{\z},H_z)}$ is the Chern
connection on $\textrm{Hom}(H_{\z},H_z)$. It is obvious that Leibniz'
rule holds for both $\deta$ and the extended $D$, with the
degree taken as the total degree in $E$, $E^{\ast}$ and $T^{\ast}(X
\times X)$. 

If $F \in \mathcal{L}^0_H$, then in analogy with Proposition \ref{ek} we have 

\begin{equation*}
D_{\textrm{Hom}(H_{\z},H_z)} \int_E F = \int_E DF.
\end{equation*}
It follows
that we also have $\debar \int_E F = \int_E \debar F$.

Let $g \in \mathcal{L}^0_H$ be such that $\nabla_{\eta} g = 0$ and
$g_{0,0}(z,z) = \textrm{Id}$.
In that case we can use $g$ as a weight just as in \eqref{thai} and get 

\begin{equation}\label{generell}
\debar  K_g = [\Delta]_{g_{0,0}} - P_g
\end{equation}
by a calculation similar to \eqref{mufti}, and then we get a Koppelman
formula by Theorem \ref{koala}.  

\begin{remark}\label{kallt}
To obtain more general formulas, one can find forms $K$ and $P$ such
that 

\begin{equation}\label{crowley}
D_{\textrm{Hom}(H_{\z},H_z)} K_g = [\Delta]_{g_{0,0}} - P_g
\end{equation}
by setting $\nabla^{\textrm{full}}_{\eta} = \deta - D$ and
checking that the corresponding equation \eqref{fredag} and Theorem
\ref{godis} are still 
valid. See for example \cite{elin} for details. This will give the
same formulas as in 
\cite{bob1}, if $H$ is the trivial line bundle. We can use weights
just as before, if we require that a weight $g$ has the property
$\nabla^{\textrm{full}}_{\eta} g = 0$ instead of $\nabla_{\eta} g = 0$. 
\end{remark}

\section{Algebraic  properties of weights}\label{viktsektion}
In this section we investigate some general constructions of weights
with the purpose of generating  
weights for a wide class of derived bundles
from two given vector bundles and weights for these. This method will
be useful later when we focus on line bundles over Grassmannians. 
 
To be more precise, we let $H$ and $H'$ be holomorphic vector 
bundles over the complex manifold $X$ and assume that $X$ fulfills the
requirements of our general setup for 
constructing Koppelman formulas, i.e., $X\times X$ admits a
holomorphic vector bundle $E$ with 
a holomorphic section defining the diagonal. Assume also that 
$g \in \Gamma(X \times X, G_{E, H})$ and 
$g' \in \Gamma(X \times X, G_{E, H'})$ are weights for 
$H$ and $H'$
respectively. 
We shall see that one can naturally define weights $g \otimes g'$ and $g \wedge g'$ (when $H=H'$), as well as 
$g^*$ for the bundles $H \otimes H', H \wedge H$ and $H^*$ respectively. This generalizes the fact, mentioned
in the introduction, that the product of weights for the trivial bundle is again a weight.

\subsection{Tensor products and exterior products of weights}
For operators $A \in H_{z} \otimes H_{\zeta}^*$ and  $B \in H_{z} \otimes (H'_{\zeta})^*$
the tensor product $A \otimes B$ defined by

\begin{equation}\label{AB}
A \otimes B (u \otimes v):=A(u) \otimes B(v), u \in H_{\zeta}, v \in H'_{\zeta} 
\end{equation}
is a linear operator in 
$\textrm{Hom}(H_{\zeta} \otimes H'_{\zeta},H_{z} \otimes H'_{z})$. 
We can therefore extend the exterior multiplication on the vector space $G_E$ 
to a linear map (which we
still denote by $\otimes$)

\begin{eqnarray*}
\otimes: (G_{E,H})_{(z,\zeta)} \otimes 
(G_{E,H'})_{(z,\zeta)}  
\rightarrow (G_{E, H \otimes H'})_{(z,\zeta)}
\end{eqnarray*}
given by

\begin{equation*}
(A \otimes \omega) \otimes (B \otimes \omega') \mapsto (A \otimes B) \otimes (\omega \wedge \omega'),
\end{equation*}
for $\omega, \omega' \in (G_E)_{(z,\zeta)}$.
This operation defines a natural fiberwise multiplication on sections.

\begin{lemma}\label{leibniz}
The operator $\nabla_{\eta}$ acts as a graded derivation with respect to the multiplication, $\otimes$, of sections, i.e.,

\begin{eqnarray*}
\nabla_{\eta} \big( (A \otimes\omega) \otimes (B \otimes \omega') \big) &=&
\nabla_{\eta}(A \otimes \omega) \otimes (B \otimes \omega') \\
& & 
+ (-1)^{deg \, \omega}(A \otimes \omega) \otimes 
\nabla_{\eta} (B \otimes \omega'),
\end{eqnarray*}
where $A$ and $B$ are local smooth sections of $H_{z} \otimes H_{\zeta}^*$ and 
$H'_{z} \otimes (H'_{\zeta})^*$
respectively, and $\omega$ and $\omega'$ are local smooth sections of $G_E$.
\end{lemma}

\begin{proof}
We first observe that 

\begin{equation*}
\nabla_{\eta}(A \otimes \omega)=-\debar A \otimes \omega+A \otimes \nabla_{\eta} \,\omega,
\end{equation*}
and likewise for $B \otimes \omega'$.
Hence,

\begin{eqnarray*}
\nabla_{\eta} \big( (A \otimes B) \otimes (\omega \wedge \omega') \big) \hspace{7cm}
\end{eqnarray*}
\begin{eqnarray*}
&=&-\debar (A \otimes B) \otimes (\omega \wedge \omega') + (A \otimes B) \otimes 
\nabla_{\eta}(\omega \wedge \omega')\\
&=&
-\debar A \otimes (B \otimes (\omega \wedge \omega')) +(A \otimes B) \otimes (\nabla_{\eta} \,\omega \wedge \omega')- \\
& &
A \otimes (\debar B \otimes \omega \wedge \omega')+ (-1)^{deg \, \omega}(A \otimes B) \otimes (\omega \wedge \nabla_{\eta} \,\omega')\\
&=&
(-\debar A \otimes \omega + A \otimes \nabla_{\eta} \,\omega) \otimes (B \otimes \omega') + \\
& & 
(A \otimes \omega ) \otimes (-\debar B \otimes \omega'+(-1)^{deg \,\omega} B\otimes \nabla_{\eta} \,\omega')\\
&=&\nabla_{\eta}(A \otimes \omega) \otimes (B \otimes \omega')+
(-1)^{deg \, \omega}(A \otimes \omega) \otimes \nabla_{\eta} (B \otimes \omega').
\end{eqnarray*}
\end{proof}

\begin{corollary}\label{tensorvikt}
Given weights $g$ and $g'$ for $H$ and $H'$ respectively, the section

\begin{equation*}
g \otimes g' \in \Gamma(X \times X, G_{E,H \otimes H'}) 
\end{equation*}
is a weight for $H \otimes H'$.
\end{corollary}

We next turn to exterior products of a vector bundle. Recall that
when $A$ and $A'$ are operators in $\mbox{Hom}(H_{\zeta}, H_z)$, $A \wedge A'$ 
is the operator in 
$\mbox{Hom}(\Lambda^2 H_{\zeta}, \Lambda^2 H_{z})$ given by 

\begin{equation*}
A \wedge A'(u \wedge u')=A(u) \wedge A'(u') - A(u') \wedge A'(u).
\end{equation*}
We can then form the exterior product

\begin{eqnarray*}
\wedge \colon (G_{E, H})_{(z,\zeta)}  \otimes (G_{E, H})_{(z,\zeta)} 
\rightarrow (G_{E, H \wedge H})_{(z,\zeta)}
\end{eqnarray*}
given by

\begin{equation*}
(A \otimes \omega) \otimes (A' \otimes \omega') \mapsto (A \wedge A') \otimes (\omega \wedge \omega').
\end{equation*}
It induces a natural exterior product on sections of $G_{E,H}$.
Using the Leibniz identity

\begin{equation*}
\debar(A \wedge A')=\debar A \wedge A' + A \wedge \debar A', 
\end{equation*}
the following lemma can be proved in the same manner as Lemma\ \ref{leibniz}.

\begin{lemma}
The operator $\nabla_{\eta}$ acts as a graded derivation with respect to the exterior multiplication of sections, i.e.,

\begin{eqnarray*}
\nabla_{\eta} \big( (A \otimes \omega) \wedge (A' \otimes \omega')\big) &=&
\nabla_{\eta}(A \otimes \omega) \wedge (A' \otimes \omega') + \\
& &
(-1)^{deg \, \omega}(A \otimes \omega) \wedge \nabla_{\eta} (A' \otimes \omega'),
\end{eqnarray*}
where $A$ and $A'$ are local smooth sections of $H_{z} \otimes H_{\zeta}^*$, and $\omega$ and 
$\omega'$ are local smooth sections of $G_E$.
\end{lemma}
In analogy with Corollary\ \ref{tensorvikt}, we have 

\begin{corollary}
Given weights $g_1$ and $g_2$ for $H$, the section
\begin{equation*}
g_1 \wedge g_2 \in  \Gamma(X \times X, G_{E,H \wedge H}) 
\end{equation*}
is a weight for $H \wedge H$.
\end{corollary}

\subsection{Dual weights}
For a local section $A \otimes \omega$ of the bundle $G_{E, H}$, 
we define the adjoint section

\begin{eqnarray*}
(A \otimes \omega)^*:=A^* \otimes \omega,
\end{eqnarray*}
where $A^*(z,\zeta)\colon H^*_{z} \rightarrow H^*_{\zeta}$ is the standard dual operator to 
$A(z,\zeta)$ given by
composing functionals with $A(z,\zeta)$. The relations

\begin{eqnarray*}
\nabla_{\eta}(A^* \otimes \omega) &=& -\debar A^* \otimes \omega +A^* \otimes \nabla_{\eta} \, \omega \\ 
&=& 
-(\debar A)^* \otimes \omega + (A \otimes \nabla_{\eta} \, \omega)^* \\
&=&
(\nabla_{\eta}(A \otimes \omega))^*
\end{eqnarray*}
prove the following lemma.

\begin{lemma}
Given a weight $g$ for the bundle $H$, the section $g^*$ is a weight for the dual bundle $H^*$.
\end{lemma}

\section{The necessary constructions on Grassmannians}\label{yuuuri}
In this section we construct the ingredients necessary to generate weighted integral formulas
on Grassmannians according to the recipe in Section \ref{yuuri}. We start by reviewing some 
elementary facts and introducing some notation. Hereafter, $X$ will denote 
the Grassmannian $Gr(k,N)$ of complex $k$-planes in $\C^N$. Just as $\CP^n$, ($=Gr(1,n+1)$),
has its tautological line bundle, $X$ has a tautological rank $k$-vector bundle, which will be 
denoted by
$H\rightarrow X$ from now on. 
We consider $H$ as a subbundle of the trivial rank $N$-bundle, $\C^N\rightarrow X$, and
the fiber of $H$ above $p\in X$ is the $k$-plane in $\C^N$ corresponding to the
point $p$. We will take the standard metric on $\C^N$ and this gives us a Hermitian metric on 
$H\subset \C^N$.
From $H$ we get a natural Hermitian line bundle $L=\det H$, which actually generates the 
Picard group; see Subsection \ref{picard}. We also get the quotient bundle, $F:=\C^N/H$, which is a 
holomorphic vector bundle of rank $N-k$. As a $C^{\infty}$-bundle,
it is isomorphic to the bundle of orthogonal complements $H^{\perp}\subset \C^N$ via the mapping
$\varphi \colon F\rightarrow H^{\perp}$ defined fiberwise by $\varphi(v+H_z)=v-\pi_{H_z} v$,
where $\pi_{H_z}$ is the orthogonal projection from $\C^N$ onto $H_z$. (If $w$ is a $\C^N$-valued
form we will, for simplicity, also write 
$\pi_{H_z} w$ for $(\pi_{H_z}\otimes \textrm{Id}) w$.)
The mapping $\varphi$ and the metric on $H^{\perp}\subset \C^N$ gives us a metric on $F$. 

Let $e=(e_1,\ldots,e_N)$ be the standard basis for $\C^N$. The point on $X$ corresponding to 
the $k$-plane $\textrm{Span}\{e_1,\ldots,e_k\}$ will be the reference point and denoted 
by $p_0$. A local holomorphic chart centered at $p_0$ can be defined as follows: Let
$z$ be a point in $\C^n:=\C^{k(N-k)}$ and organize $z$ as an $(N-k)\times k$-matrix, i.e.,

\begin{equation*}
z=\left(
\begin{array}{ccc}
z_{11} & \cdots & z_{1k} \\
\vdots &         & \vdots \\
z_{N-k,1} & \cdots & z_{N-k,k}
\end{array}
\right) \in \C^{n}.
\end{equation*}
Associate to $z$ the point on $X$ corresponding to the $k$-plane spanned by the columns
of the $N\times k$-matrix

\begin{equation}\label{korv}
\left(
\begin{array}{c}
I \\
z
\end{array}
\right), \,\,\, I=I_{k\times k},
\end{equation}
with respect to the basis $e$.
This actually gives us an injective map from $\C^{n}$ onto a dense subset $U\subset X$.
We also get natural local holomorphic frames for the bundles $H$, $L$, and $F$ over this chart.
For $j=1,\ldots,k$, let $\Hram_j(z)$ be the $j$th column of \eqref{korv}, i.e., 
$\Hram_j(z)=e_j+\sum_{i=1}^{N-k}z_{ij}e_{k+i}$. Then $\Hram_1,\ldots,\Hram_k$ are 
$k$ pointwise linearly independent holomorphic sections of $H$ over $U$. A natural holomorphic 
frame for $L$ is thus $\mathfrak{l}=\Hram_1 \wedge \cdots \wedge \Hram_k$. 
Also, for $1\leq j \leq N-k$, let $\Fram_j(z)$ be the equivalence class defined by $e_{k+j}$  
in $F=\C^N/H$, in the fiber over $z$. Then $(\Fram_1,\ldots,\Fram_{N-k})$ is a local 
holomorphic frame for $F$ over $U$. The projection $\C^N\rightarrow F$, expressed in the
$e$-basis for $\C^N$ and the frame $\Fram$ for $F$, can then be written as the 
$(N-k)\times N$-matrix

\begin{equation}\label{kvotproj}
\left(
\begin{array}{cc}
-z & I
\end{array}
\right), \,\,\, I=I_{(N-k)\times (N-k)}.
\end{equation}
For reference we note some more explicit expressions:
As a mapping $\C^N_e\rightarrow \C^N_e$ expressed in the $e$-basis we have

\begin{equation*}
\pi_{H}=\left(
\begin{array}{c}
I \\
z
\end{array}
\right)
(I+z^*z)^{-1}
\left(
\begin{array}{cc}
I & z^*
\end{array}
\right)
\end{equation*}
and as a mapping $\C^N_e\rightarrow H_{\Hram}$,

\begin{equation*}
\pi_{H}=(I+z^*z)^{-1}
\left(
\begin{array}{cc}
I & z^*
\end{array}
\right).
\end{equation*}
The mapping $\varphi \colon F_{\Fram}\rightarrow \C^N_e$ looks like

\begin{equation*}
\varphi=
\left(
\begin{array}{c}
-(I+z^*z)^{-1}z^* \\
I-z(I+z^*z)^{-1}z^*
\end{array}
\right).
\end{equation*}
We have defined the metric, $\langle \cdot , \cdot\rangle_F$, on $F$ via $\varphi$
so the Hermitian metric-matrix, $h_F$, expressed in the frame $\Fram$ satisfies
$(h_F)_{i,j}=\langle \varphi(\Fram_i),\varphi(\Fram_j)\rangle_{\C^N}$, (with the convention
that $\langle v,w\rangle_F=v^th_F\bar{w}$). 
Using the explicit expression for $\varphi$, a computation then gives

\begin{equation*}
h_F^t(z)=(I+zz^*)^{-1},
\end{equation*}
and so the Chern curvature-matrix of $F$ is

\begin{equation*}
\Theta_F=\debar (\bar{h}_F^{-1}\partial \bar{h}_F) = \partial \debar \log (I+zz^*),
\end{equation*} 
where the last expression should be interpreted in the functional calculus sense. 
For the bundle $H$ we get

\begin{equation*}
h_H^t=I+z^*z, \,\,\, \textrm{and}\,\,\,\, \Theta_H=\partial \debar \log (I+z^*z)^{-1},
\end{equation*}
expressed in the frame $\Hram$.

\subsection{The bundle $E$ and the section $\eta$}
We will construct a holomorphic vector bundle $E\rightarrow X_z \times X_{\z}$ of rank
$n$ ($=k(N-k)$)
and a global holomorphic section $\eta$ of it defining the diagonal. As in Section \ref{yuuri}, 
we let $H_z$ and $H_{\z}$
denote the pull-back of the tautological bundle under the projections 
$X_z\times X_{\z}\rightarrow X_z$ and $X_z\times X_{\z}\rightarrow X_{\z}$ respectively and we
define $F_z$ similarly. However, for convenience we will occasionally abuse this notation
and also write, e.g., $H_z$ for the fiber of the bundle $H_z\rightarrow X_z\times X_{\z}$ above a point $(z,\z)$.
This ambiguity is (partly) justified since one can identify fibers of $H_z \rightarrow X_z\times X_{\z}$ above points 
$(z,\z)$ for any $\z$. This means also that, e.g, $\{\Hram_j (z)\}$ is a local holomorphic frame for 
$H_z\rightarrow X_z\times X_{\z}$ over $U_z\times X_{\z}$.

The bundle $E$ is simply $E=F_z\otimes H_{\z}^*$ and then
$\Eram_{ij}:=\Fram_i(z)\otimes \Hram_j^*(\z)$, $1\leq i\leq N-k$, $1\leq j \leq k$, is a 
holomorphic frame for $E$ over $U\times U \subset X\times X$. To define $\eta$ we start
with a vector $v\in H_{\z}$ and via $H_{\z}\subset \C^N_{\z}\cong \C^N_z$ we can identify $v$ 
with a vector $\tilde{v}\in \C^N_z$. We then let $\eta(v)$ be the projection of $\tilde{v}$ on
$F_z=\C^N_z/H_z$. 

\begin{proposition}
The section $\eta$ of $E$ is holomorphic and defines the diagonal in $X\times X$.
\end{proposition}

\begin{proof}
It is clear that $\eta(v)$ vanishes if and only if $v$ belongs to the fiber
above a point in the diagonal $\Delta\subset X\times X$.
Hence, $\eta$ is a global section of $\Hom(H_{\z},F_z)\cong E$ and vanishes precisely on $\Delta$.
In the coordinates and frames described above, $\eta$ has the form

\begin{equation*}
\eta=\z-z.
\end{equation*} 
In fact, if $v=\sum_1^k v_j\Hram_j(\z)\in H_{\z}$ then $\eta(v)$ is the image in 
$F_z$ of $\sum_{1}^kv_je_j+\sum_{i=1}^{N-k}\sum_{j=1}^k\z_{ij}v_je_{k+i}$. By \eqref{kvotproj}
this is equal to $\sum_{i=1}^{N-k}\sum_{j=1}^k(\z_{ij}-z_{ij})v_je_{k+i}$.
We thus see that $\eta$ is holomorphic and vanishes to the first order on $\Delta$.
\end{proof}

\subsection{Bundles and weights}\label{bil}
The bundle $L=\det H$ actually generates the Picard group of holomorphic line bundles; 
cf.\ Section 5.3, and \cite{snow}. 
We will construct weights for the line bundles $L^r:=L^{\otimes r}\rightarrow X$,
and for the vector bundle $H\rightarrow X$. We start 
by defining two fundamental sections $\gamma_0$ and $\gamma_1$ of 
$\Hom(H_{\z},H_z)$ and $\Hom(H_{\z},H_z)\otimes E^*\wedge T^*_{0,1}(X \times X)$
respectively. For $v\in H_{\z}$ we first identify $v$ 
with the vector $\tilde{v}$ in the trivial bundle $\C^N_z\rightarrow X_z\times X_{\z}$
via $H_{\z}\subset \C^N_{\z} \cong \C^N_z$. We then put $\gamma_0(v)=\pi_{H_z}\tilde{v}$. 
In the $\Hram$-frames described above, $\gamma_0$ is simply the $k\times k$-matrix

\begin{equation}\label{gamma0}
\gamma_0=(I+z^*z)^{-1}(I+z^*\z).
\end{equation}
It is a little bit more complicated to describe $\gamma_1$: Let $\xi$ and $v$ be (germs of) 
smooth sections of $E$ and $H_{\z}$ respectively. Since $E=F_z\otimes H_{\z}^*$,
$\xi(v)$ defines naturally a smooth section of $F_z$ and hence,
$\varphi(\xi(v))$ is a smooth section of $H_z^{\perp}\subset \C^N_z$.
We then put $-\gamma_1(\xi\otimes v)=\pi_{H_z}(\debar \, \varphi(\xi(v)))$, which
is a smooth section of $H_z\otimes T^*_{0,1}(X\times X)$. We check that $\gamma_1$ so
defined actually is tensorial. Let $f$ be (a germ of) a smooth function. We then get

\begin{eqnarray*}
\gamma_1(f\xi\otimes v) &=& -\pi_{H_z}\,\big( \debar \, \varphi(f \xi(v))\big) \\
&=& -\pi_{H_z} \,\big( \varphi(\xi(v))\otimes \debar f+ f\debar \varphi(\xi(v))\big) \\
&=& -\pi_{H_z}\big( \varphi(\xi(v))\big) \otimes \debar f + f\gamma_1(\xi\otimes v).
\end{eqnarray*}
But $\pi_{H_z}\big( \varphi(\xi(v))\big)=0$ since $\varphi(\xi(v)) \in H_z^{\perp}$, and so
$\gamma_1(f\xi\otimes v)=f\gamma_1(\xi \otimes v)$. (One could also note that 
$\gamma_1(\xi\otimes v)=-[\pi_{H_z},\debar]\varphi (\xi(v))$, where $[\pi_{H_z},\debar]$ is the 
commutator.) Hence, $\gamma_1$ defines a section of
$H_z\otimes T^*_{0,1}(X_z\times X_{\z})\otimes E^* \otimes H^*_{\z}\cong 
\Hom(H_{\z},H_z)\otimes E^* \wedge T^*_{0,1}(X\times X)$. 
A computation in the local coordinates shows that 

\begin{equation}\label{gamma11}
\gamma_1=\sum_{i,j=1}^k\Hram_i(z)\otimes \Hram_j^{*}(\z)\otimes M_{ij},
\end{equation}
where $M$ is the $k\times k$-matrix of $E^{*}$-valued $(0,1)$-forms

\begin{equation}\label{gamma12}
M=\debar\big((I+z^{*}z)^{-1}z^{*}\big)\wedge \Eram^{*}.
\end{equation}
Here, $\Eram^{*}$ is the matrix with entries $(\Eram_{ij})^*$.

\begin{proposition}\label{viktpropp}
The section $G:=\gamma_0 + \gamma_1\in \mathcal{L}^0_H$, (cf.\ \eqref{vikter}),
is a weight for the tautological bundle $H$.
\end{proposition}

\begin{proof}
We need to check that $\gamma_0(z,z)=\textrm{Id}$ and that $\nabla_{\eta} G=0$.
The first equality is obvious from the definition. For the second one we have to verify the two
equations $\debar \gamma_0=\deta \gamma_1$ and $\debar \gamma_1=0$.
Let $v$ be a germ of a holomorphic section of $H_{\z}$. Via $H_{\z}\subset \C^N_{\z}\cong \C^N_z$
we may view $v$ as a holomorphic section of $\C^N_z$ and then we can write 

\begin{eqnarray*}
(\deta \gamma_1)(v) &=& -\pi_{H_z}\big(\debar(\varphi(\eta(v)))\big)=
-\pi_{H_z} \big(\debar(\pi_{H_z^{\perp}}v)\big) \\
&=&
-\pi_{H_z} \big(\debar(v-\pi_{H_z}v)\big) =
\debar_{H_z}(\pi_{H_z}v) \\
&=& \debar_{H_z}(\gamma_0(v)).
\end{eqnarray*}
Hence, $\debar_{H_z}(\gamma_0(v))=\deta \gamma_1(v)$ for any germ of holomorphic section
$v$ of $H_{\z}$. It follows that $\debar \gamma_0=\deta \gamma_1$.
Now, let $\xi$ be a germ of a holomorphic section of $E$. Then $\xi(v)$ is a germ of a holomorphic
section of $F_z$. One can (locally) lift $\xi(v)$ to a germ of a holomorphic section, $\widetilde{\xi(v)}$,
of $\C^N$ that projects to $\xi(v)$. We then get

\begin{eqnarray*}
\debar_{H_z} \gamma_1(\xi\otimes v) &=& \debar_{H_z} \big(\pi_{H_z}\debar(\varphi(\xi(v)))\big) 
= \debar_{H_z} \big(\pi_{H_z}\debar(\widetilde{\xi(v)}-\pi_{H_z}\widetilde{\xi(v)})\big) \\
&=& 
-\debar_{H_z} \big(\pi_{H_z}\debar(\pi_{H_z}\widetilde{\xi(v)})\big) 
=-\debar^2_{H_z}(\pi_{H_z}\widetilde{\xi(v)}) \\
&=& 0.
\end{eqnarray*}
Hence, $\debar_{H_z} \gamma_1(\xi\otimes v)=0$ for any holomorphic $\xi$ and $v$, and this finishes the proof.
\end{proof}

By the algebraic properties of weights established in Section \ref{viktsektion} we now get that 
$g:=G\wedge \cdots \wedge G$ (the exterior product of $G$ with itself $k$ times) is a weight 
for $L$.
It is easy 
to check that 

\begin{equation*}
g_{0,0}=\gamma_0\wedge\cdots \wedge \gamma_0 =
k!\frac{\det (I+z^*\z)}{\det (I+z^*z)}
\end{equation*}
in the frame $\mathfrak{l}$ for $L$.
Weights for positive powers of $L$ are then obtained by taking powers of $g$. By the results at 
the end of Section \ref{viktsektion} we can also get weights for $H^*$ and $L^{-r}=(L^*)^{\otimes r}$
from $G$. If one wants to construct weights for $H^*$ geometrically, as we have done in this section,
it is easier to take $F_{\z}\otimes H_{z}^*$ as the bundle $E$.
However, our Koppelman formulas have an inherent duality and this gives us weighted 
formulas for forms with values in $H^*$ and $L^{-r}$ from the weighted formulas for $H$ and $L^r$.

\section{Representation-theoretic interpretations}\label{yuuuuri}
In this section we describe $X$ and its line bundles in terms of group actions and representations. 
The purpose of this is threefold. First of all, this point of view gives an easy description of the Picard 
group of $X$. Secondly, and more importantly, we prove that the weights we have constructed earlier
will all be invariant under a certain group action; a property which will turn out be highly useful
in the last section with applications to Bergman kernels.
Finally, in this setup, we can fairly easily prove that the restriction of the bundle $E$ to the diagonal
is equivalent to the holomorphic cotangent bundle $T^*_{1,0}$ of $X$.

\subsection{The Grassmannian as a homogeneous space}
The linear action of the group $GL(N,\mathbb{C})$ on $\mathbb{C}^N$ induces an action as 
holomorphic automorphisms of $X$, 
and this action is clearly transitive.
Hence, we can describe $X$ as a homogeneous space $X \cong GL(N,\mathbb{C})/P$, where

\begin{eqnarray*}
P:=\left\{\left.\left( \begin{array}{cc} A & B \\0 & D \end{array}\right) \right| \det A \det D \neq 0\right\}
\end{eqnarray*}
is the stabilizer of $p_0$. 
One can also restrict the action to the subgroup $SL(N, \mathbb{C})$ and still 
have a transitive group action; this time
exhibiting $X$ as the homogeneous space $SL(N,\mathbb{C})/P'$, where 

\begin{eqnarray*}
P':=\left\{\left.\left( \begin{array}{cc} A & B \\0 & D \end{array}\right) \right| \det A \det D =1\right\}
\end{eqnarray*}
is the stabilizer of $p_0$ in $SL(N,\mathbb{C})$. The reason that we mention this realization is that some of the 
results we refer to later hold only for quotients of semisimple Lie groups.
A third realization is given by restricting the $GL(N,\mathbb{C})$-action to the unitary group $U(N)$. 
The stabilizer of $p_0$ in this subgroup is

\begin{eqnarray*}
\left\{\left.\left( \begin{array}{cc} A & 0 \\0 & D \end{array}\right) \right| A \in U(k), D \in U(N-k)\right\} 
\cong U(k) \times U(N-k),
\end{eqnarray*} 
and hence we have a third description of $X$ as the quotient space $U(N)/(U(k) \times U(N-k))$.

\subsection{The bundles $H$, $F$,  and $E$}
We recall that a vector bundle $\mathcal{V} \rightarrow X$ is said to be \emph{homogeneous} 
under a group $G$ if $G$ acts on it by 
bundle automorphisms in such a way that the corresponding action on $X$ is transitive. As a 
consequence, the stabilizer, $G_{p_0}$, of
$p_0$ in $G$ acts linearly on the fiber $\mathcal{V}_{p_0}$, i.e., $\mathcal{V}_{p_0}$ 
carries a representation, $\tau$, of $G_{p_0}$.
The vector bundle $\mathcal{V}$ can then be reconstructed from the representation $\tau$ 
as the set of equivalence classes

\begin{eqnarray*}
G \times_{G_{p_0}} \mathcal{V}_{p_0}:=G \times \mathcal{V}_{p_0} / \sim,
\end{eqnarray*}    
where the equivalence relation $\sim$ is defined as $(g,v) \sim (g',v')$ if 
and only if $(g',v')=(gx^{-1},\tau(x)v)$ for some $x$ in  $G_{p_0}$.
The $G$-action is then given by $[(g',v)] \stackrel{g}{\mapsto} [(gg',v)]$, 
where the brackets denote the equivalence classes of the 
respective pairs. The holomorphic vector bundles are those associated with 
holomorphic representations, $\tau$, of $G_{p_0}$, i.e., 
$\tau: G_{p_0} \rightarrow \mbox{End}(\mathcal{V}_{p_0})$ is a holomorphic group homomorphism.

Suppose now that $H \subset G$ is a closed subgroup of $G$ which also acts 
transitively on $X$. Then we can describe $X$ as a quotient 
$H/(H \cap G_{p_0})$ and form the $H$-homogeneous vector bundle 
$\mathcal{V}^H:=H \times_{H \cap G_{p_0}} \mathcal{V}_{p_0}$.
This latter bundle is in fact equivalent to the former one via the bundle mapping

\begin{eqnarray*}
\Psi^G_H: H \times_{H \cap G_{p_0}} \mathcal{V}_{p_0} & \rightarrow &  G \times_{G_{p_0}} \mathcal{V}_{p_0},\\
\left[(h,v)\right]_{H} &\mapsto& \left[(h,v)\right]_G,
\end{eqnarray*}
where the brackets denote the respective equivalence classes.

For our purposes, this means that we can choose to view $GL(N,\mathbb{C})$-homogeneous 
vector bundles as $SL(N,\mathbb{C})$-homogeneous ones 
without 
any loss of information as long as the corresponding representations 
of $P'$ are restrictions of $P$-representations.
Moreover, since $P$ is the complexification of $U(k) \times U(N-k)$ 
(i.e., $U(k) \times U(N-k)$ is a totally real submanifold of $P$), 
a holomorphic representation of $P$ is uniquely determined by its restriction 
to $U(k) \times U(N-k)$. Hence we can also view the
vector bundle as only $U(N)$-homogeneous. 

The group $GL(N, \mathbb{C})$ acts naturally on the trivial bundle $X \times \mathbb{C}^N$ by 
$(p,v) \stackrel{g}{\mapsto} (g(p),gv)$. The tautological bundle $H$ is invariant 
under this action, and is therefore a 
$GL(N,\mathbb{C})$-homogeneous vector bundle. We let 
$\tau: P \rightarrow \mbox{End}(\mathbb{C}^k)$ denote the 
corresponding representation of $P$ on $H_{p_0}\cong\mathbb{C}^k$, namely

\begin{eqnarray*}
\tau\left(\begin{array}{cc}A & B \\ 0 & D \end{array} \right)v=Av, \,v \in \mathbb{C}^k. 
\end{eqnarray*}
Since the subbundle $H$ of $\mathbb{C}^N$ is $GL(N,\mathbb{C})$-invariant, there is a well-defined action 
on the quotient bundle $F=\mathbb{C}^N/H$; i.e., $F$ is also a homogeneous bundle.
We can identify the fiber $F_{p_0}$ with $\mathbb{C}^{N-k}$, and we let $\rho$ 
denote the corresponding $P$-representation given by

\begin{eqnarray*}
\rho\left(\begin{array}{cc}A & B \\ 0 & D \end{array} \right)v=Dv, \,v \in \mathbb{C}^{N-k}.
\end{eqnarray*}
The bundle $E \rightarrow X \times X$ is homogeneous under the product group 
$GL(N,\mathbb{C}) \times GL(N,\mathbb{C})$, and the representation
of $P \times P$ on the fiber 
$(F_{z} \otimes H_{\zeta}^*)_{(p_0,p_0)} \cong \mbox{Hom}(\mathbb{C}^k, \mathbb{C}^{N-k})$ is the tensor product
representation $\rho \otimes \tau^*$ given by

\begin{eqnarray*}
\rho \otimes \tau^*(g_{z},g_{\zeta})Z=D_{z}ZA_{\zeta}^{-1},\quad 
g_{\zeta}&=&\left(\begin{array}{cc}A_{\zeta} & B_{\zeta} \\ 0 & D_{\zeta} \end{array} \right),\\
g_{z}&=&\left(\begin{array}{cc}A_{z} & B_{z} \\ 0 & D_{z} \end{array} \right),\\
Z &\in& M_{N-k,k}(\mathbb{C}).
\end{eqnarray*} 

The trivial bundle $\mathbb{C}^N$ is equipped with the standard Euclidean metric 
which is $U(N)$-invariant; and the tautological 
bundle $H$ inherits this metric, thus admitting an isometric action of $U(N)$. 
Moreover, we recall that the quotient bundle $F$ is 
smoothly equivalent to the orthogonal complement, $H^{\perp}$, to the tautological bundle.
It should be pointed out that $H^{\perp}$ is not a holomorphic vector bundle, whereas $F$ is. 
Since the metric on $F$ is induced from that on $H^{\perp}$,
the $U(N)$-action on $F$ is also isometric. Moreover, the bundle $E$ is equipped with a tensor product metric, and
therefore the Cartesian product $U(N) \times U(N)$ acts isometrically on $E$.

The Chern connections and curvatures of the three bundles $H, F$, and $E$ are invariant 
under the respective group actions
since they are associated with invariant metrics. We recall that the invariance of a curvature, 
$\Theta_{\mathcal{V}}$, of 
a holomorphic homogeneous vector bundle $\mathcal{V}$ means the invariance as 
a section of the bundle $\mbox{End}(\mathcal{V}) \otimes T^*_{1,1}$ with respect to the 
natural action on sections of this 
bundle. Concretely, this means that

\begin{eqnarray*}
\Theta_{\mathcal{V}}(gp)(u,v) w=g \Theta_{\mathcal{V}}(p)(dg^{-1}(gp)u,dg^{-1}(gp)v)g^{-1}w,\\
u \in T^*_{(1,0),gp},\quad
v \in T^*_{(0,1),gp},\quad
w \in \mathcal{V}_{gp}.
\end{eqnarray*}
In particular, it follows that the curvature is determined by its value at a fixed reference point. We shall return
to the Chern curvature of $E$ below, and give an explicit formula for it at the point $p_0$. First, however, we shall
undertake a closer study of the restriction of $E$ to the diagonal.

The action of the group $U(N)$ on $X$ defines a fibration $q:U(N) \rightarrow X$ given by 
$q(g)=g(p_0)$ which is $U(N)$-equivariant
with respect to left multiplication, $L_g:x \mapsto gx$, on the group itself, 
and the action on $X$, i.e., 
$q(gx)=g(q(x))$ holds for $g, x \in U(N)$.
Moreover, the right action $R_l:x \mapsto xl^{-1}$ of the subgroup $U(k) \times U(N-k)$ 
on $U(N)$ preserves each fiber $q^{-1}(p)$ for 
$p \in X$, and yields a diffeomorphism $U(k) \times U(N-k) \cong q^{-1}(p)$. This equips $U(N)$ with
the structure of a principal $U(k) \times U(N-k)$-bundle over $X$. Since the right action 
of $U(k) \times U(N-k)$ commutes with
left multiplication, the group $U(N)$ acts equivariantly with respect to the action of 
$U(k) \times U(N-k)$.
Moreover, the embedding of $U(N)$ into $M_N(\mathbb{C})$ induces an Riemannian structure 
on $U(N)$ by restriction of the trace inner product
$(A,B) \mapsto \mbox{tr}(AB^*)$, and the left multiplication is isometric with respect 
to this inner product. 
For any $g \in U(N)$ with $q(g)=p$, we have an orthogonal decomposition

\begin{eqnarray}
T_g(U(N))=T_g(q^{-1}(p)) \oplus T_g(q^{-1}(p))^{\perp}, \label{dekomp}
\end{eqnarray}
and this decomposition is invariant under left multiplication.
The restriction of the differential of $q$ to the orthogonal complement 
$T_g(q^{-1}(p))^{\perp}$ yields an isomorphism

\begin{eqnarray*}
dq(g)|_{T_g(q^{-1}(p))^{\perp}}:T_g(q^{-1}(p))^{\perp} \rightarrow T_{q(g)}(X).
\end{eqnarray*}
For any $p \in X$ we thus have a family of subspaces parametrized by the set 
$q^{-1}(p)$ to which the tangent space at $p$ is isomorphic. 
We therefore define an equivalence relation on the tangent bundle $T(U(N))$ by

\begin{eqnarray}
(g,v) \sim (g',v') \quad \mbox{iff} \quad (g',v')=(R_l(g),dR_l(g)v),  \label{ekvi}
\end{eqnarray}
for some $l \in U(k) \times U(N-k)$. 
By the isometry of the left multiplication, the orthogonal complement bundle $\cup_p T(q^{-1}(p))^{\perp}$ is
a $U(N)$-homogeneous vector bundle. Moreover, for any vector in this subbundle, 
the whole equivalence class lies in the subbundle since
also the right action is isometric. 
It follows that $S:=\cup_p T(q^{-1}(p))^{\perp}/\sim$ is a well-defined $U(N)$-homogeneous vector bundle over $X$.
Clearly, $S$ is equivalent to the tangent bundle $T(X)$, and thus it inherits a complex structure.

\begin{proposition}\label{Epropp}
The restriction of $E$ to the diagonal $\Delta(X \times X)$ is equivalent to the 
holomorphic cotangent bundle $T^*_{1,0}(X)$.
\end{proposition}

\begin{proof}
We prove that $E^*$ is equivalent to $S$. Since $S$ is $U(N)$-homogeneous, it is 
uniquely determined by the corresponding representation
of $U(k) \times U(N-k)$ on the fiber $S_{p_0}$. For the identity element $e \in U(N)$, 
the tangent space $T_e(U(N))$ is
isomorphic to the Lie algebra

\begin{eqnarray*}
\mathfrak{u}(N):=\{ X \in M_N(\mathbb{C})|X^*=-X\},
\end{eqnarray*}
and the subspaces in the decomposition \eqref{dekomp} are 
explicitly given by

\begin{eqnarray}
T_e(q^{-1}(p_0))&=&\left\{ 
\left.\left( \begin{array}{cc} Y & 0\\ 0 & Z\end{array}\right)\right| Y^*=-Y, Z^*=-Z  \right\},\\
T_e(q^{-1}(p_0))^{\perp}&=&\left\{ \left.\left( \begin{array}{cc} 0 & B\\ -B^* & 0\end{array}\right)\right| 
B \in M_{k,N-k}(\mathbb{C})  \right\}. \label{tangentrum}
\end{eqnarray}
For $v= \left( \begin{array}{cc} 0 & B\\ -B^* & 0\end{array}\right) \in T_e(q^{-1}(p_0))^{\perp}$, and 
$l =\left( \begin{array}{cc} A & 0\\ 0 & D\end{array}\right) \in U(k) \times U(N-k)$, 

\begin{eqnarray*}
dL_l(e)v&=&\left( 
\begin{array}{cc} A & 0\\ 0 & D\end{array}\right)\left( \begin{array}{cc} 0 & B\\ -B^* & 0\end{array}\right)\\
&=&\left( \begin{array}{cc} 0 & AB \\ -DB^* & 0\end{array}\right).
\end{eqnarray*}
We can represent the equivalence class of this tangent vector by a tangent vector at the identity, namely by

\begin{eqnarray*}
dR_{l}(l)dL_l(e)v&=&\left( \begin{array}{cc} 0 & AB \\ -DB^* & 0\end{array}\right)
\left( \begin{array}{cc} A^{-1} & 0\\ 0 & D^{-1}\end{array}\right)\\
&=&\left( \begin{array}{cc} 0 & ABD^{-1}\\ -(ABD^{-1})^* & 0 \end{array}\right).
\end{eqnarray*}
Hence, we can identify the representation of $U(k) \times U(N-k)$ on $S_{p_0}$ 
with the representation on $M_{k,N}(\mathbb{C})$ given 
by $B \mapsto ABD^{-1}$, i.e., with the representation $\tau \otimes \rho^*$ on 
$\mbox{Hom}(\mathbb{C}^{N-k}, \mathbb{C}^k)$, which
is precisely the $U(k) \times U(N-k)$ representation associated to the restriction of $E^*$ to the diagonal.
\end{proof}

\begin{remark}
In \cite{bob1}, Berndtsson proves that any appropriate bundle $E\rightarrow X\times X$ 
has to coincide with the holomorphic 
cotangent bundle on the diagonal. In the case of $\CP^n$, an independent proof of Proposition \ref{Epropp} 
can be found in the book \cite{Demailly} by Demailly; Proposition 15.7 in Chapter V. 
\end{remark}

By the identification $T_{p_0}(X)$ with the subspace $T_e(q^{-1}(p_0))^{\perp}$ in 
\eqref{tangentrum}, we have an explicit realization of 
its complexification

\begin{eqnarray*}
T_{p_0}(X)^{\mathbb{C}} \cong \left\{ \left.\left( \begin{array}{cc} 0 & B\\ C & 0\end{array}\right)\right| 
B \in M_{k,N-k}(\mathbb{C}), C \in M_{N-k,k}(\mathbb{C})  \right\}.
\end{eqnarray*}
Consider now the element  
$\left( \begin{array}{cc} i \frac{N-k}{N}I_k & 0\\ 0 & -i\frac{k}{N}I_{N-k}\\ \end{array}\right) \in 
T_e(q^{-1}(p_0)) \cong \mathfrak{u}(k) \times \mathfrak{u}(N-k)$. 
Its adjoint action determines the complex structure, $J_{p_0}$, at $p_0$ by

\begin{eqnarray*}
J_{p_0}\left( \begin{array}{cc} 0 & B\\ -B^* & 0\end{array}\right)&:=&
\left[\left( \begin{array}{cc} i \frac{N-k}{N}I_k & 0\\ 0 & -i\frac{k}{N}I_{N-k}\\ \end{array}\right),
\left( \begin{array}{cc} 0 & B\\ -B^* & 0\end{array}\right)\right]\\
&=&\left( \begin{array}{cc} 0 & iB\\ -(iB)^* & 0\end{array}\right).
\end{eqnarray*}
The splitting of $T_{p_0}^{\mathbb{C}}$ into the the $\pm i$-eigenspaces is given by

\begin{eqnarray*}
T_{(1,0),p_0}(X) &\cong& \left\{ \left.\left( \begin{array}{cc} 0 & Y\\ 0 & 0\end{array}\right)\right| 
Y \in M_{k,N-k}(\mathbb{C}) \right\},\\
T_{(0,1),p_0}(X) &\cong& \left\{ \left.\left( \begin{array}{cc} 0 & 0\\ Z & 0\end{array}\right)\right| 
Z \in M_{N-k,k}(\mathbb{C})  \right\}.
\end{eqnarray*} 
We recall that the curvature $\Theta_E$ at the point $p_0$ is given by the formula

\begin{eqnarray}
\Theta_E(p_0)(Y,Z)W=(\rho \otimes \tau^*)'([Y,Z])(W), \label{kurv}
\end{eqnarray}
where $(\rho \otimes \tau^*)'$ denotes the differentiated representation of 
the Lie algebra $\mathfrak{u}(k) \times \mathfrak{u}(N-k)$ given
by  $(\rho \otimes \tau^*)'(X):=\frac{d}{dt} (\rho \otimes \tau^*)(\exp tX)|_{t=0}$. 
The explicit expression for \eqref{kurv} is

\begin{eqnarray*}
\Theta_E(p_0)(Y,Z)W&=&(\rho \otimes \tau^*)'\left( \begin{array}{cc} YZ & 0\\ 0 & ZY\end{array}\right)(W)\\
&=&ZYW-WYZ, \quad W \in M_{N-k,k}(\mathbb{C}).
\end{eqnarray*}

\subsection{Invariance of weights}\label{invofweights}
In this section we study a natural action of $U(N)$ on sections of the bundles 
$\mbox{Hom}(L^r_{\zeta},L^r_z) \otimes G_E$, and prove
that the corresponding weights are invariant under that action.

Recall that for an action of a group, $G$, on a vector bundle $\mathcal{V} \rightarrow M$, 
a natural action is induced on the space of sections by 

\begin{equation}\label{snittverkan}
(gs)(z):=gs(g^{-1}z),
\end{equation} 
where the second action on the right hand side refers to the action on the total space of the bundle.
The bundles  $\mbox{Hom}(L^r_{\zeta},L^r_z) \otimes G_E$ are equipped with the natural $U(N) \times U(N)$ 
actions given as tensor (and exterior) products of the actions described in the previous 
section and their duals. In what follows, we will consider the action of $U(N)$ 
(embedded as the diagonal subgroup of $U(N) \times U(N)$) given by restriction. 
The actions on the respective total spaces are the obvious ones, and 
we will therefore use the simple notation from \eqref{snittverkan} for such an action.

We let $g^r:=g^{\otimes r}$ for $r \geq 0$ and $g^r:=(g*)^{\otimes r}$ for $ r \leq 0$ 
denote the weight for the line bundle $L^r$.

\begin{proposition}\label{invariantvikt}
The weight $g^r$ is a $U(N)$-invariant section of the vector bundle $\mbox{Hom}(L^r_{\zeta},L^r_z) \otimes G_E$.
\end{proposition}

\begin{proof}
It clearly suffices to prove that the section $G=\gamma_0+\gamma_1$ is an invariant section of 
$\mbox{Hom}(H_{\zeta},H_{\zeta})$;
and for this, we prove that $\gamma_0$ and $\gamma_1$ are invariant separately.
We now fix an orthonormal basis, $\{h_1, \ldots, h_k\}$, for $H_z$. For any $u \in H_{\zeta}$ and $l \in U(N)$, we have

\begin{eqnarray*}
(l\gamma_0)(u)=l\gamma_0(l^{-1}u)
&=&l \sum_{i=1}^k \langle l^{-1}u,l^{-1}h_i \rangle l^{-1}h_i \\
&=& \sum_{i=1}^k \langle u,h_i \rangle h_i\\
&=&\gamma_0(u),
\end{eqnarray*}
which shows the invariance of $\gamma_0$.
We now consider $\gamma_1$, and therefore choose a local section $f$ of $F$ near the point $z \in X$. 
Then, we have

\begin{eqnarray*}
(l\gamma_1)(f \otimes u)&=&-l(\pi_{H_{l^{-1}z}}(\debar \varphi(l^{-1}f)) \otimes l^{-1}u)\\
&=&-l(\pi_{H_{l^{-1}z}}(l^{-1}\debar \varphi(f)) \otimes l^{-1}u)\\
&=&-\pi_{H_z}(\debar \varphi(f) \otimes u)\\
&=&\gamma_1(f \otimes u),
\end{eqnarray*}
where the third equality is completely analogous to the invariance of $\gamma_0$.
This concludes the proof.
\end{proof}

We now turn our attention to the form $P_{g^r}$ defined in \eqref{thai} again. 

\begin{corollary}
The form $P_{g^r}$ is $U(N)$-invariant.
\end{corollary}

\begin{proof}
First of all, an argument similar to the proof of Proposition\ \ref{invariantvikt} shows that the 
section $\eta$ is $U(N)$-invariant.
Secondly, the Chern connection $D_E$ on $E$ commutes with the $U(N)$-action, and hence $D\eta$ is also $U(N)$-invariant.
The curvature $\Theta$ is even $U(N) \times U(N)$-invariant; and hence it follows that the form
$g \wedge \left(\frac{D \eta}{2 \pi i} + \frac{i \tilde{\Theta}}{2 \pi}\right)_n$ is $U(N)$-invariant.
We now claim that the operator $\int_E$ is $U(N)$-equivariant. Indeed, the identity section $I \in \mbox{End}(E)$ is 
obviously $U(N)$-invariant, and so is therefore also the section $\tilde{I}_n$ defined in connection 
with Definition \ref{E-integral}. 
Hence, $\int_E$ is an equivariant operator, and this also finishes the proof.   
\end{proof}

The canonical splitting $T^*(X \times X) \cong T^*_z(X) \oplus T^*_{\zeta}(X)$ of the cotangent bundle 
of $X \otimes X$ is $U(N) \times U(N)$-invariant, 
and hence $(P_{g^r})$ can be decomposed as 

\begin{equation}
(P_{g^r})=\sum_{\stackrel{\large{p'+p''=n}}{\small{q'+q''=n}}}(P_{g^r})_{p',p'',q',q''}, \label{fyrindex}
\end{equation}
where $(P_{g^r})_{p',p'',q',q''}$ is a section of 
$\mbox{Hom}(H_{\zeta},H_z) \otimes \Lambda^{p',q'}(T^{\mathbb{C}}_z)^* \wedge \Lambda^{p'',q''}(T^{\mathbb{C}}_{\zeta})^*$,
i.e., it is of bidegree $(p',q')$ in the $z$-variable, and of bidegree $(p'',q'')$ in the 
$\zeta$-variable according to the splitting. By the invariance of the splitting, we also have 

\begin{corollary}\label{invcor}
The terms $(P_{g^r})_{p',p'',q',q''}$ in the decomposition \eqref{fyrindex} are $U(N)$-invariant.
\end{corollary}

Only the term $(P_{g^r})_{n,0,n,0}$ which has bidegree $(n,n)$ in the $z$-variable will 
contribute to the integral in the Koppelman formula.
Later we will examine this term more closely. 

\begin{corollary}
The current $K_{g^r}$ in \eqref{thai} is $U(N)$-invariant.
\end{corollary}

\begin{proof}
It clearly suffices to prove that $u$ in \eqref{nord} is $U(N)$-invariant; and since the group action commutes with
the $\debar$-operator and exterior powers, it only remains to prove the invariance of $\sigma$.
Note that $\sigma$ can be described by the equation

\begin{equation*}
\sigma (v)=\frac{\langle v,\eta \rangle_E}{|\eta|_E^2}, \quad v \in E.
\end{equation*}
The invariance of $\sigma$ now follows immediately from the invariance of $\eta$ and from
the fact that the action of $U(N)$ preserves the metric.
\end{proof}

\subsection{Line bundles on $X$}\label{picard}
In this subsection we recapitulate how the Picard group of $X$ can be described in terms of 
holomorphic characters. All of this is classical theory
and well-known, even though the results in their explicit form can be hard to find in the literature.
The reason for including it in the paper is rather to give an  overview for readers who are not 
familiar with representation theory of Lie groups.

Suppose now that $\mathcal{L} \rightarrow X$ is a $SL(N,\mathbb{C})$-homogenous holomorphic 
line bundle. The corresponding $P'$-representation
then amounts to a holomorphic character $\chi_{\mathcal{L}}:P' \rightarrow \mathbb{C}^*$. 
Moreover, it is well-known that all holomorphic line
bundles over $X$ are in fact $SL(N,\mathbb{C})$-homogeneous (cf. \cite{snow}), and hence 
the Picard group $H^1(X,\mathcal{O}^*)$ is isomorphic to 
the multiplicative group of holomorphic characters of $P'$.\\

Suppose now first that $\chi \colon P \rightarrow \mathbb{C}^*$ is a holomorphic character. 
(This is no restriction, as we
shall later see that all holomorphic characters of $P'$ are restrictions of $P$-characters.) 
It is well-known that it is 
then uniquely determined by
its restriction to the Levi-subgroup $GL(k) \times GL(N-k)$ realized as 

\begin{eqnarray*}
\left\{\left.\left(\begin{array}{cc}
A & 0\\
0 & D\\
\end{array}
\right)\right| \det A \det D \neq 0 \right\}.
\end{eqnarray*} 
By restricting to the respective factors, we can uniquely express $\chi$ as a product 
$\chi=\chi_1\chi_2$, where $\chi_1$ and $\chi_2$
are characters of $GL(k,\mathbb{C})$ and $GL(N-k,\mathbb{C})$ respectively. Let 
$\chi_1': \mathfrak{gl}(k,\mathbb{C}) \rightarrow \mathbb{C}$
denote the differential at the identity of $\chi_1$. Then $\chi_1'$ annihilates 
the commutator ideal in the decomposition 

\begin{eqnarray*}
\mathfrak{gl}(k,\mathbb{C}) = \mathfrak{Z}(\mathfrak{gl}(k,\mathbb{C})) \oplus \left[\mathfrak{gl}(k,\mathbb{C}), 
\mathfrak{gl}(k,\mathbb{C})\right]
\end{eqnarray*}
of $\mathfrak{gl}(k,\mathbb{C})$ as the direct sum of the center and the commutator. 
More specifically, we have the identity

\begin{eqnarray*}
 \left[\mathfrak{gl}(k,\mathbb{C}), \mathfrak{gl}(k,\mathbb{C})\right]=\mathfrak{sl}(k,\mathbb{C}),
\end{eqnarray*}
from which it follows that the normal subgroup $SL(k,\mathbb{C})$ lies in the kernel of the 
character $\chi_1$. Hence, $\chi_1$ descends to
a character, $\widetilde{\chi_1}$, of the quotient group $GL(k,\mathbb{C})/SL(k, \mathbb{C})$, 
yielding the commuting diagram

\begin{eqnarray*}
\xymatrix{GL(k,\mathbb{C}) \ar[r]^{\chi_1} \ar@{->>}[d]& \mathbb{C}^*\\
GL(k,\mathbb{C})/SL(k,\mathbb{C}) \ar[ur]_{\widetilde{\chi_1}}}.
\end{eqnarray*}
Moreover, the quotient group is isomorphic to $\mathbb{C}^*$ via the mapping 
$gSL(k,\mathbb{C}) \mapsto \det g$, and hence we have the diagram

\begin{eqnarray*}
\xymatrix{GL(k,\mathbb{C}) \ar[r]^{\chi_1} \ar@{->>}[d]& \mathbb{C}^*\\
GL(k,\mathbb{C})/SL(k,\mathbb{C}) \ar[ur]_{\widetilde{\chi_1}} \ar@{^{(}->>}[r]  & \mathbb{C}^*, \ar[u]}
\end{eqnarray*}
which allows us to identify $\widetilde{\chi_1}$ with a holomorphic character 
$\mathbb{C}^* \rightarrow \mathbb{C}^*$.
The latter ones are easily described. Indeed, by holomorphy, any such character 
is uniquely determined by its restriction to the totally real subgroup
$S^1 \subset \mathbb{C}^*$, on which it gives a character $S^1 \rightarrow S^1$. 
Hence, it is of the form $\zeta \mapsto \zeta^m$, for some integer $m$. 
The analogous result holds of course for $\chi_2$. Summing up, we have thus found that 

\begin{eqnarray*}
\chi\left( \begin{array}{cc}
A & B\\
0 & D\\
\end{array}\right)
=\det A^m \det D^n, 
\end{eqnarray*}
for some $m, n \in \mathbb{Z}$.

The line bundle corresponding to the choice $m=1, n=0$ is the determinant of the tautological vector bundle. 
To study the line bundle corresponding to the parameters $m=0, n=1$, we consider it as 
a $SL(N, \mathbb{C})$-homogeneous line bundle,
which amounts to restricting the corresponding character to the subgroup $P'$ of $P$. 
We let $\chi'$ denote
the differential at the identity of this character. The Lie algebra 
$\mathfrak{p}'$ admits a decomposition

\begin{eqnarray*}
\mathfrak{p}'=\mathfrak{Z}(\mathfrak{p}') \oplus \left[\mathfrak{p}', \mathfrak{p}'\right]
\end{eqnarray*}
as the direct sum of its center and its commutator ideal. These two ideals are given by

\begin{eqnarray*}
\mathfrak{Z}(\mathfrak{p}')&=&\left\{ \left.\left( \begin{array}{cc} c (N-k) I_k & 0\\ 0 & -ckI_{N-k}\\ \end{array}\right)\right| 
c \in \mathbb{C}\right\},\\
 \left[\mathfrak{p}', \mathfrak{p}'\right]&=&\left\{ \left.\left( \begin{array}{cc} A  & B\\ 0 & D\\ \end{array}\right)\right| 
\mbox{tr} A=\mbox{tr} D=0\right\}.
\end{eqnarray*}
On the group level, we have the commutator subgroup 

\begin{eqnarray*}
[P',P']=\left\{ \left. \left(\begin{array}{cc} 
A & B\\
0 & D\\
\end{array}\right) \right| \det A=\det D=1\right\},
\end{eqnarray*}
and the quotient group $P'/[P',P']$ has complex dimension one. In fact, an isomorphism 
$\Phi:P'/[P',P'] \rightarrow \mathbb{C}^*$ is given by

\begin{eqnarray*}
\Phi(g\,[P',P'])=\det A,
\end{eqnarray*}
for $g=\left(\begin{array}{cc}
A & B\\
0 & D
\end{array}
\right)$.

If $\mu: P' \rightarrow \mathbb{C}^*$ is a holomorphic character, it factors through 
the projection onto the quotient group just as above, yielding
a holomorphic character $\tilde{\mu}: P'/[P',P'] \rightarrow \mathbb{C}^*$. Using the 
isomorphism $\Phi$ above, we obtain the commuting diagram

\begin{eqnarray*}
\xymatrix{P' \ar[r]^{\mu} \ar@{->>}[d]& \mathbb{C}^*\\
P'/[P',P'] \ar[ur]_{\widetilde{\mu}} \ar@{^{(}->>}[r]  & \mathbb{C}^*. \ar[u]}
\end{eqnarray*}
From this, we conclude that $\mu \left( \begin{array}{cc}
A & B\\
0 & D \\
\end{array}\right)=\det A^j$, for some $j \in \mathbb{Z}$. In particular, 
it follows that $\mu$ can naturally be extended to
a holomorphic character $P \rightarrow \mathbb{C}^*$. Moreover, the dual 
bundle to the determinant of the tautological vector bundle
corresponds to the $P'$-character $\left( \begin{array}{cc}
A & B\\
0 & D \\
\end{array}\right) \mapsto \det A^{-1}=\det D$, which can be extended to the $P$-character  $\left( \begin{array}{cc}
A & B\\
0 & D \\
\end{array}\right) \mapsto \det D$.
It is easy to see that the $GL(N,\mathbb{C})$-homogeneous line bundle associated with this holomorphic character 
is isomorphic to the determinant of
the quotient bundle $F=\mathbb{C}^N/H$.

\subsection{The Bott-Borel-Weil theorem}
In this subsection we briefly describe some group representations associated with homogeneous vector bundles.

Suppose now that $G$ is a complex Lie group acting transitively and holomorphically on a complex manifold $M$, so that 
we can write $M \cong G/T$ for some closed subgroup $T \subseteq G$. Let $\mathcal{V} \rightarrow M$ be a $G$-homogeneous
holomorphic vector bundle. Recall that the action of $G$ on $\mathcal{V}$ induces the action on smooth sections
given by \eqref{snittverkan}. 
Since $G$ acts holomorphically on $M$, there is a natural action on $\mathcal{V}$-valued $(p,q)$-forms 
(by taking the pullback composed with inversion). Moreover, the action commutes with the $\overline{\partial}$-operator on 
$\mathcal{V}$, from which it follows that the action preserves closed forms and exact form; thus inducing an
action on the Dolbeault cohomology groups $H^{p,q}(M,\mathcal{V})$. In the case when $G$ is a complexification
of some semisimple compact Lie group, $G_{\mathbb{R}}$, the Bott-Borel-Weil theorem 
(cf.\ \cite{penrose}, Theorem. 5.0.1) 
gives a realization of all irreducible representations of  $G_{\mathbb{R}}$ as $H^{0,q}(M,\mathcal{L})$ 
for some homogeneous 
line bundle, $\mathcal{L}$, over $M$, and also states the vanishing of the other sheaf 
cohomology groups associated with $\mathcal{L}$.
We shall see examples of it in the context of the vanishing theorems of the next section.

\section{Applications}\label{yuuuuuri}
\subsection{Vanishing theorems}

We would like to find vanishing theorems for the bundles $L^r$ and
$L^{-r}$ over $X$ by means of the Koppelman
formula. This will yield explicit solutions to the
$\debar$-equation in the cohomology groups which are trivial. 

Let $D$ in Theorem \ref{koala} be the whole of $X$, and let
$\phi(\z)$ be a 
$\debar$-closed form of bidegree $(p,q)$ taking values in $L_\z^r$, with
$r>0$. The only obstruction to solving the $\debar$-equation is then the 
term $\int_\z \phi(\z) \wedge P_{g^r}(\z,z)$. We have 

\begin{eqnarray} \label{tuktoyaktuk}
& P_{g^r} = & \int_{E} g^r \wedge \left(\frac{D \eta}{2 \pi
    i} + \frac{i \tilde{\Theta}}{2 \pi}\right)_n = \\
& = & \int_{E} \sum_{j = 1}^{\min(kr,n)} C_{j} (g^r)_{j,j} \wedge
(D \eta)^j \wedge (\tilde{\Theta}_{\z} + \tilde{\Theta}_z)^{n-j} \nonumber
\end{eqnarray}
where $(g^r)_{j,j}$ is the term in $g^r$ which has bidegree $(0,j)$
and takes values in $\Lambda^j E^{\ast}$. Note that all the
differentials in $g$ are in the $z$ variable; this is because
$\debar_\z$ commutes with $\pi_{H_z}$. 

\begin{theorem} \label{manzanita}
The cohomology groups $H^{p,q}(X, L^r)$ are trivial in the following cases:
\begin{itemize}
\item[a)] $p \neq q$ and $r = 0$.

\item[b)] $p > q$ and $r > 0$.

\item[c)] $p < q$, $rk < q - p$, and $r > 0$.

\item[d)] $p < q$ and $r < 0$.

\item[e)] $p > q$, $rk < p - q$, and $r < 0$.
\end{itemize}
\end{theorem}

\begin{proof}
\emph{a)} If $r=0$ we do not need a weight, and in that
case 

\begin{equation*}
P = \int_E \left( \frac{i\tilde{\Theta}}{2 \pi} \right)_n = c_n(E),
\end{equation*}
or the $n$:th Chern form of $E$. It is obvious that $P$ consists of terms
with bidegree $(k,k)$ in $z$ and $(n-k,n-k)$ in $\z$, and thus $\int
\phi \wedge P = 0$ if $\phi$ has bidegree $(p,q)$ with $p \neq q$. \\
\\
\noindent
\emph{b)} Since the only source of antiholomorphic
differentials in $\z$ is $\tilde{\Theta}_\z$, which is a $(1,1)$-form,
we can never get more $d \bar{\z}_i$:s than $d \z_i$:s. This means that
$\int_\z \phi(\z) \wedge P_{g^r} = 0$ if $\phi$ has bidegree $(p,q)$
where $p > q$ (since then $P_{g^r}$ would need to have bidegree
$(n-p,n-q)$ in $\z$ with $n-q > n-p$). \\
\\
\noindent
\emph{c)} If $\phi(\z)$ has bidegree $(p,q)$, then $P_{g^r}$
needs to have bidegree $(p,q)$ in $z$. We can take at most $p$ of the
$\tilde{\Theta}_{z}$:s. We will then need at least $q - p$ more $d
\bar{z}_i$:s,  and these have to come from the factor $g^r$. But $g^r$
has maximal bidegree $(0,rk)$, so if $rk < q-p$ the obstruction will
vanish. \\  
\\
\noindent
\emph{d)} By duality, if we have a $(p,q)$-form
$\phi$ taking values in $L^r$ with $r<0$, the obstruction is given by
$\int_z \phi(z) \wedge P_{g^{-r}}(\z,z)$. This is zero unless there is a term in
$P_{g^{-r}}$ of bidegree $(p,q)$ in $\z$. By the same argument as in the proof of
b), the obstruction vanishes if $q > p$. \\  
\\
\noindent
\emph{e)} If $\phi(z)$ has bidegree $(p,q)$, then
$P_{g^{-r}}$
needs to have bidegree $(n-p,n-q)$ in $z$, where $n-q > n-p$. The
rest follows as in the proof of c). 
\end{proof}

\begin{remark}
In $\CP^n$, we can get rid of the obstruction in more cases, either by
proving that $P_{g^r}$ is 
$\debar_\z$-exact (since then Stokes' theorem can be applied), or by
proving that it is $\debar_z$-exact 
(since then $\int_\z \phi \wedge P_{g^r}$ will be $\debar_z$-exact as
well). See \cite{elin} for details. 
\end{remark}

Part $d)$ of the above theorem is the special case of the Bott-Borel-Weil theorem for the 
parabolic quotient $GL(N,\mathbb{C})/P$. 
For $r=-1$, all vanishing theorems were proved by le Potier in \cite{lePotier}. He also proved vanishing
theorems for exterior and symmetric powers of the tautological bundle and its dual. In \cite{snow0}, Snow
gives an algorithm for computing all Dolbeault cohomology groups for all line bundles over Grassmannians.
Implementing the algorithm in a computer, Snow obtains various vanishing theorems including ours.
It is worth noting that both le Potier and Snow obtain their results by reduction to the Bott-Borel-Weil theorem.

\subsection{Bergman kernels}
We will see that the projection part, $P_{g^r}$, of our Koppelman formula for $L^r$ basically is 
the Bergman kernel associated with the space of holomorphic sections of $L^{-r}$. 
We begin by examining $P_{g^r}$. Recall that 

\begin{equation*}
g^r=\big((\gamma_0+\gamma_1)^k\big)^{\otimes r}=
\big(\sum_{j=0}^k\binom{k}{j}\gamma_0^{k-j}\wedge \gamma_1^j\big)^{\otimes r}
=:(\gamma_0^k)^{\otimes r}+\tilde{g}^r,
\end{equation*}
where $\gamma_0^k$ of course is the $k$th exterior power of $\gamma_0$, is our weight for $L^r$. 
The projection kernel in our Koppelman formula for $L^r$ is thus

\begin{eqnarray*}
P_{g^r}&=&\int_E g^r\wedge \left( \frac{D\eta}{2\pi i}+\frac{i\tilde{\Theta}_E}{2\pi} \right)_n \\
&=& (\gamma_0^k)^{\otimes r}\otimes \int_E \left( \frac{D\eta}{2\pi i}+\frac{i\tilde{\Theta}_E}{2\pi}\right)_n
+\int_E \tilde{g}^r\wedge \left( \frac{D\eta}{2\pi i}+\frac{i\tilde{\Theta}_E}{2\pi}\right)_n \\
&=:& P_{g^r}^0+\tilde{P}_{g^r}.
\end{eqnarray*}
Let $\PV_{g^r}$ and $\PV^0_{g^r}$ be the parts of $P_{g^r}$ and $P^0_{g^r}$ respectivly, which have bidegree
$(n,n)$ in the $z$-variables. 
Let us examine $\PV_{g^r}$ and $\PV^0_{g^r}$ more closely on the set $Z:=\{p_0\}\times X_{\z}$.
In our local coordinates and frames over $U_z\times U_{\z}$ we have by 
\eqref{gamma0} that $\gamma_0=(I+z^*z)^{-1}(I+z^*\z)$.
On $Z$ intersected with $\{p_0\}\times U_{\z}$, denoted $Z'$ below, we thus have $\gamma_0=I$ expressed in our frames.
According to \eqref{gamma11} and \eqref{gamma12}, 
we see that, as a matrix in our frames for $H_z$ and $H_{\z}$, $\gamma_1=dz^*\wedge \Eram^*$
on $Z'$. Moreover, a straightforward computation shows that
the part of $D\eta$, which does not contain any differentials in the $\z$-variables,
equals $-\sum_{i,j}dz_{ij}\wedge \Eram_{ij}$ on $Z'$. Also, the part of $\tilde{\Theta}_E$,
which does not contain any differentials in the $\z$-variables, is 
$\widetilde{\Theta_{F_z}\otimes \textrm{Id}_{H_{\z}^*}}$. We thus see that the building blocks for 
$\PV_{g^r}$ and $\PV^0_{g^r}$ are independent of $\z$ on $Z'$ 
when expressed in our frames. Since both $\PV_{g^r}$ and $\PV^0_{g^r}$
take values in a line bundle we must have $\PV_{g^r}=C\PV^0_{g^r}$ on $Z'$. But $Z'$ is dense in $Z$
and so this equality holds on $Z$ by continuity.
Now, by Corollary \ref{invcor} in Subsection \ref{invofweights}, it follows that both
$\PV_{g^r}$ and $\PV^0_{g^r}$ are invariant under the diagonal group in $U(N)\times U(N)$ and
since $Z$ intersects each orbit under this group
we can conclude that $\PV_{g^r}=C\PV^0_{g^r}$ on all of $X\times X$.

Given any holomorphic section $f$ of $L^{-r}$, $r>0$, and any vector $v_p$ in the fiber of
$L^r$ above an arbitrary point $p$, our Koppelman formula now gives

\begin{equation}\label{Beq1}
f(p).v_p = \int_{X_z}\PV_{g^r}(z,p)\wedge v_p\wedge f(z)=
C\int_{X_z}\PV^0_{g^r}(z,p)\wedge v_p\wedge f(z).
\end{equation}
It is easy to compute $\PV^0_{g^r}$ explicitly, and one gets

\begin{equation*}
\PV^0_{g^r}= \big(\frac{i}{2\pi}\big)^n (\gamma_0^k)^{\otimes r}\otimes
\int_E \left( \widetilde{\Theta_{F_z}\otimes \textrm{Id}_{H^{*}_{\z}}}\right)_n=
\big(\frac{i}{2\pi}\big)^n (\gamma_0^k)^{\otimes r}\otimes c_{N-k}(\Theta_{F_z})^k.
\end{equation*}
Moreover, $\Theta_{F_z}$ is the $U(N)$-invariant curvature of $F_z$, so it follows that 
$c_{N-k}(\Theta_{F_z})^k$ is a $U(N)$-invariant $(n,n)$-form and hence equal to a constant times the
invariant volume form $dV$. We have thus obtained 

\begin{equation}\label{Beq2}
f(p).v_p=C\int_{X_z}f(z).(\gamma_0^k)^{\otimes r}v_p dV(z) 
\end{equation}
for any holomorphic section $f$ of $L^{-r}$. Modulo a multiplicative constant, one
also has that $dV=(c_1(L))^n$, and then the above formula assumes the following form
expressed in the frames and coordinates discussed above.

\begin{equation*}
f(\z)=C\int_{\C^n}f(z)\frac{\det (I+z^{*}\z)^r}{\det (I+z^{*}z)^r} 
\big((\partial \debar \log \det (I+z^*z))\big)^n.
\end{equation*}

We will now describe what will be the Bergman kernel. 
Let $\rho^r \colon L_z^r \rightarrow L_z^{-r}$ be the antilinear identification induced by the metric,
i.e., $\rho^r(v)=\langle \cdot , v \rangle_{L^r_z}$, and define
$K_r(z,\z)\colon L_{\z}^r\rightarrow L_z^{-r}$ by $K_r(z,\z)=\rho^r \circ (\gamma_0^k)^{\otimes r}$.
Then one easily checks that $K_r(z,\z)$ is a fiberwise antilinear map, which depends antiholomorphically
on $\z$. To show that it actually depends holomorphically on $z$ we consider the adjoint operator 
$K_r(z,\z)^{*}\colon L^r_z \rightarrow L_{\z}^{-r}$ and the operator 
$K_r(\z,z)\colon L^r_z \rightarrow L_{\z}^{-r}$. We know that the latter operator depends 
antiholomorphically on $z$. Note also that since $K_r(z,\z)$ is fiberwise antilinear, the adjoint 
should be defined by $(K_r(z,\z)^*u). v=\overline{u.(K_r(z,\z)v)}$ for 
$u\in L^{r}_z$ and $v\in L^r_{\z}$. It is then straightforward to check that 
$K_r(z,\z)^{*}=K_r(\z,z)$, and so $K_r(z,\z)^{*}$ must depend antiholomorphically on $z$. It follows that
$K_r(z,\z)$ depends holomorphically on $z$. In particular, for any non-zero vector $v\in L_p^r$, the mapping
$z\mapsto K_r(z,p)v$ defines a global non-zero 
holomorphic section of $L^{-r}$.
In fact, these sections generate $H^0(X,L^{-r})$ as we now show. Consider the Bergman space
$A^2_r$ defined as $H^0(X,L^{-r})$ equipped with the norm

\begin{equation*}
\|f\|_{A^2_r}^2:=\int_{X}\|f\|^2_{L^{-r}}dV,\,\,\, f\in
H^0(X,L^{-r}).
\end{equation*}

We claim that, modulo a multiplicative constant,  
$K_r(z,\z)$ is the Bergman kernel for $A^2_r$, i.e., that $K_r(z,\z)$
is the fiberwise antilinear map 
$L^r_{\z}\rightarrow L^{-r}_z$, which depends holomorphically
on $z$ and antiholomorphically on $\z$, and has the property that for any $f\in A^2_r$ and any vector
$v\in L^r_{\z}$ (in the fiber above $\z$) one has
\begin{equation*}
f(\z).v=\langle f,K_r(\cdot,\z)v\rangle_{A^2_r}=
\int_{X}\langle f(z),K_r(z,\z)v\rangle_{L^{-r}_z}dV(z).
\end{equation*}
It only remains to verify this last property. But this reproducing property 
follows directly from \eqref{Beq2} after noting the following equality, which basically is the 
definition of $K_r(z,\z)$:

\begin{equation*}
u. ((\gamma_0^k)^{\otimes r}v)=\langle u, K_r(z,\z)v \rangle_{L^{-r}_z}, \,\,\,
\textrm{for all}\,\, u\in L^{-r}_z, \,\, \textrm{and all} \,\, v\in L^r_{\z}.
\end{equation*}

\begin{remark}
In the case of $\CP^n$ it is not too hard to compute $\mathcal{P}_{g^r}$ directly from its definition. 
For instance,
one can first verify in local or homogeneous coordinates that the part of 
$\gamma_1\wedge D\eta$ which contains no $d\z$ or $d\bar{\z}$ is equal to
$-\gamma_0\otimes \widetilde{\Theta_{F_z}\otimes \textrm{Id}}_{\mathcal{O}(1)_{\z}}$,
cf.\ Proposition 4.1 and the weight $\alpha$ in \cite{elin}. Then, a straightforward 
computation shows that $\mathcal{P}_{g^r}$ is equal to

\begin{equation*}
\binom{n+r}{n}\big(\frac{i}{2\pi}\big)^n\gamma_0^r\otimes \det (\Theta_{F_z}).
\end{equation*}
\end{remark}


\begin{thebibliography}{99}

\bibitem{MA1} \textsc{M. Andersson:} Integral representation with weights
I. \textit{Math. Ann.}, \textbf{326(1)} (2003), 1--18. 

\bibitem{penrose} \textsc{R. J. Baston and M. G. Eastwood:}
The {P}enrose transform: Its interaction with representation theory.
\textit{Oxford Science Publications}, Clarendon Press, Oxford, 1989.
    
\bibitem{BY2} \textsc{C. Berenstein, R. Gay, A. Vidras and A. Yger:} Residue
    currents and Bezout identities. \textit{Progress in Mathematics 114},
Birkhuser Verlag, Basel, 1993.


\bibitem{bob1} \textsc{B. Berndtsson:} Cauchy-Leray forms and vector bundles.
\textit{Ann. Sci. cole Norm. Sup. (4)}, \textbf{24(3)} (1991), 319--337.

\bibitem{BE2} \textsc{B. Berndtsson:} Integral formulas on projective 
space and the Radon transform of Gindikin-Henkin-Polyakov. \textit{Publ. 
Mat.}, \textbf{32(1)} (1988), 7--41.

\bibitem{AB} \textsc{B. Berndtsson and M. Andersson:} Henkin-Ramirez
    formulas with weight factors.  \textit{Ann. Inst. Fourier (Grenoble)}, \textbf{32(3)}
  (1982), v--vi, 91--110. 

\bibitem{Demailly} \textsc{J.-P. Demailly:} Complex Analytic and Differential Geometry. 
Online book, available at: http://www-fourier.ujf-grenoble.fr/$\sim$demailly/.

\bibitem{elin} \textsc{E. G\"{o}tmark:} Weighted integral formulas on manifolds. 
\textit{Ark. Mat.}, to appear. Available at ArXiv: math.CV/0611082.

\bibitem{HL} \textsc{G. M. Henkin and J. Leiterer:} Global integral formulas
for solving the $\bar \partial $-equation on Stein manifolds,
\textit{Ann. Polon. Math.}, \textbf{39} (1981), 93--116. 


\bibitem{hua} \textsc{L. K. Hua:}
Harmonic analysis of functions of several complex variables in the classical domains.
\textit{Translations of Mathematical Monographs},
\textbf{6}, American Mathematical Society, Providence, R.I., 1979.


\bibitem{K-W} \textsc{A. W. Knapp and N. R. Wallach:}
Szeg\"o kernels associated with discrete series.
\textit{Invent. Math.}, {\bf 34(3)} (1976), 163--200.

\bibitem{lang} \textsc{R. P. Langlands:}
Dimension of spaces of automorphic forms.
\textit{Algebraic Groups and Discontinuous Subgroups (Proc. Sympos.
              Pure Math., Boulder, Colo., 1965)},
253--257, Amer. Math. Soc., Providence, R.I., 1966.
  
\bibitem{lePotier} \textsc{J. le Potier:} Cohomologie de la grassmannienne  valeurs dans les puissances
extrieures et symtriques du fibr universel.
\textit{Math. Ann.}, \textbf{226(3)} (1977), 257--270.
    	

\bibitem{QU} \textsc{D. Quillen:} Superconnections and the Chern character.
\textit{Topology}, \textbf{24(1)}, (1985), 89--95.


\bibitem{satake} \textsc{I. Satake:} Algebraic structures of symmetric domains.
\textit{Kan\^o Memorial Lectures},
\textbf{4}, Iwanami Shoten, Tokyo, 1980.


\bibitem{schmid-L^2} \textsc{W. Schmid:} $L\sp{2}$-cohomology and the discrete series.
\textit{Ann. of Math. (2)},
{\bf 103(2)} (1976), 375--394.
     

\bibitem{snow0} \textsc{D. M. Snow:} Cohomology of twisted holomorphic forms on Grassmann
manifolds and quadric hypersurfaces.
\textit{Math. Ann.}, \textbf{276(1)} (1986), 159--176. 	

\bibitem{snow} \textsc{D. M. Snow:} Homogeneous vector bundles.
\textit{Group actions and invariant theory (Montreal, PQ, 1988)},
CMS Conf. Proc., \textbf{10}, 193--205, Amer. Math. Soc.,
Providence, RI, 1989.

\bibitem{wong1} \textsc{H.-W. Wong:} Dolbeault cohomological realization of Zuckerman 
modules associated with finite rank representations.
\textit{J. Funct. Anal.},
\textbf{129} (1995), 428--454..

\bibitem{wong2} \textsc{H.-W. Wong:}
Cohomological induction in various categories and the maximal globalization conjecture.
\textit{Duke Math. J.},
\textbf{96(1)} (1999), 1--27.

\bibitem{genkai} \textsc{G. Zhang:}
Berezin transform on compact Hermitian symmetric spaces.
\textit{Manuscripta Math.}, \textbf{97(3)} (1998), 371--388.    

\end{thebibliography}
\end{document}